\documentclass[11pt,reqno]{amsart}
\usepackage{}
\usepackage{a4wide}
\numberwithin{equation}{section}
\usepackage{mathrsfs}
\usepackage{amsfonts}
\usepackage{amsmath}
\usepackage{stmaryrd}
\usepackage{amssymb}
\usepackage{amsthm}
\usepackage{mathrsfs}
\usepackage{url}
\usepackage{amsfonts}
\usepackage{amscd}
\usepackage{indentfirst}
\usepackage{enumerate}
\usepackage{amsmath,amsfonts,amssymb,amsthm}
\usepackage{amsmath,amssymb,amsthm,amscd}
\usepackage{graphicx,mathrsfs}
 \usepackage{appendix}

\usepackage[colorlinks,linkcolor=blue, citecolor=blue]{hyperref}
\numberwithin{equation}{section}

\usepackage[numbers,sort&compress]{natbib}

\usepackage{esint}
\usepackage{graphicx}
\usepackage[dvipsnames]{xcolor}
\newcommand{\R}{\mathbb{R}}

\newcommand\e\varepsilon
\newtheorem{Thm}{Theorem}[section]
\newtheorem{Lem}{Lemma}[section]
\newtheorem{Prop}{Proposition}[section]

\newtheorem{Rem}{Remark}[section]

\begin{document}
\title[Pohozaev identities for weak solutions of Grushin type p-sub-Laplacian equation]
{Pohozaev identities for weak solutions of Grushin type p-sub-Laplacian equation via domain variations}

\author[Y. Wei and X. Zhou]{Yawei Wei and Xiaodong Zhou}

\address[Yawei Wei]{School of Mathematical Sciences and LPMC, Nankai University, Tianjin 300071, China}
\email{weiyawei@nankai.edu.cn}

\address[Xiaodong Zhou]{School of Mathematical Sciences, Nankai University, Tianjin 300071, China}
\email{1120210030@mail.nankai.edu.cn}

\thanks{Acknowledgements: This work is supported by the NSFC under the grands 12271269 and the Fundamental Research Funds for the Central Universities.}





\begin{abstract} In this paper, we study Pohozaev identities for weak solutions of degenerate elliptic equations involving Grushin type p-sub-Laplacian under only $C^1$-regularity assumption. By using domain variations, we obtain the local Pohozaev identities of translating type and scaling type. As an application, a global Pohozaev identity of scaling type in $\R^{N+l}$ is also derived.

\end{abstract}

\maketitle

\section{Introduction}

\setcounter{equation}{0}

For any $p>1$, we denote Grushin type p-sub-Laplacian as follows
\begin{equation}\label{x5}
  \Delta_\gamma^pu=div_\gamma(|\nabla_\gamma u|^{p-2}\nabla_\gamma u).
\end{equation}
In which, $\nabla_\gamma$ is the Grushin gradient defined by
\begin{equation}\label{x2}
  \nabla_\gamma u=(\nabla_xu,|x|^\gamma\nabla_yu),
\end{equation}
with
\begin{equation}\label{x0}
  z=(x,y)=(x_1,\cdots,x_N,y_1,\cdots,y_l)\in\R^{N+l},
\end{equation}
and $N+l\geq3$. $\gamma\geq0$ is a real number. $div_\gamma$ is the Grushin divergence given by
\begin{equation}\label{x4}
  div_\gamma(\overrightarrow{P}(z))=\sum_{i=1}^N\frac{\partial p_i(z)}{\partial x_i}+|x|^\gamma\sum_{j=1}^l\frac{\partial q_j(z)}{\partial y_j},
\end{equation}
for vector field
\begin{equation}\label{x3}
  \overrightarrow{P}(z)=(p_1(z),\cdots,p_N(z),q_1(z),\cdots,q_l(z)).
\end{equation}
Especially,
\begin{equation}\label{x1}
  \Delta_\gamma u=div_\gamma(\nabla_\gamma u)=\Delta_xu+|x|^{2\gamma}\Delta_yu,
\end{equation}
is the well-known Grushin operator.

In this paper, we consider the following Dirichlet problem
\begin{equation}\label{x6}
\begin{cases}
 -\Delta_\gamma^pu=f(z,u),\;\;&~\mbox{in}~\Omega,\\[2mm]
u=0, ~ \;\;&~\mbox{on}~\partial\Omega,
\end{cases}
\end{equation}
where $\Omega$ is a bounded open domain in $\R^{N+l}$. The nonlinearity $f:\R^{N+l}\times\R^1\rightarrow\R^1$ is a continuous function. Let
\begin{equation}\label{zzxxyyy}
  F(z,u)=\int_0^uf(z,s)ds.
\end{equation}
A function $u\in W_{\lambda,0}^{1,p}(\Omega)$ is a weak solution of \eqref{x6} if
\begin{equation}\label{x7}
  \int_\Omega|\nabla_\gamma u|^{p-2}\nabla_\gamma u\cdot\nabla_\gamma\varphi dz=\int_\Omega f(z,u)\varphi dz,\;\;\forall\varphi\in W_{\lambda,0}^{1,p}(\Omega).
\end{equation}
In which, $W^{1,p}_{\lambda,0}(\Omega)$ is the closure of $C_0^\infty(\Omega)$ in $W^{1,p}_\lambda(\Omega)$ with respect to the norm $\|\cdot\|_{W^{1,p}_\lambda(\Omega)}$. And for $p\geq1$, we denote
\begin{equation}\label{x8}
  W^{1,p}_\lambda(\Omega)=\left\{u\in L^p(\Omega)\;|\;\frac{\partial u}{\partial x_i}, |x|^\gamma\frac{\partial u}{\partial y_j}\in L^p(\Omega), i=1,\cdots,N, j=1,\cdots,l.\right\},
\end{equation}
and
\begin{equation}\label{x9}
  \|u\|_{W^{1,p}_\lambda(\Omega)}=\big(\int_{\Omega}|\nabla_\gamma u|^p+|u|^pdz\big)^\frac{1}{p}.
\end{equation}
We define
\begin{equation}\label{x420}
  g(z)=f(z,u(z)).
\end{equation}
Then the energy functional corresponding to \eqref{x6} is:
\begin{equation}\label{x10}
  I(u)=\frac{1}{p}\int_\Omega|\nabla_\gamma u|^pdz-\int_\Omega g(z)udz,\;\;u\in W^{1,p}_{\lambda,0}(\Omega).
\end{equation}
For the definition \eqref{x7} of weak solution and the definition \eqref{x8} of weighted Sobolev space, we refer to \cite{zengjia1}. And the Grushin type p-sub-Laplacian is also involved in \cite{zengjia4}.

Pohozaev identity plays a very important role in proving the existence and nonexistence results for
the nonlinear elliptic partial differential equations. Pohozaev identities for several typical elliptic equations on bounded domains or unbounded domains were studied thoroughly in \cite{Cao14} by Cao and his collaborators. This is an overview article. After establishing this kind of identity, mathematicians mentioned in this review proved many existence and nonexistence results of the solutions by using Pohozaev identities. For example,
\begin{equation}\label{i4}
\begin{cases}
-\Delta u=\lambda f(u),~x\in\Omega,\\[2mm]
u|_{\partial\Omega}=0,
\end{cases}
\end{equation}
where $\Omega\subset\R^N$ is an any bounded open domain of $\R^N$ ($N\geq3$), $f(u):\R^1\rightarrow\R^1$ is a continuous function satisfies $f(0)=0$. If $u\in C^2(\Omega)\cap C^1(\overline{\Omega})$ is the solution of \eqref{i4}, then the following Pohozaev identity for scaling holds:
\begin{equation}\label{i5}
  \frac{1}{2}\int_{\partial\Omega}(x\cdot\nu)|\nabla u|^2dS=\frac{2-N}{N}\lambda\int_\Omega uf(u)dx+N\lambda\int_\Omega F(u)dx,
\end{equation}
where $F(u)=\displaystyle\int_0^uf(t)dt$ and $\nu=\nu(x)$ is the unit outer normal vector of the point $x\in\partial\Omega$.
Identity \eqref{i5} established by Pohozaev in \cite{Pohozaev15} initially is for the sake of prove \eqref{i4} only has trivial solution under some technical conditions. However, the proof of such identities in \cite{Pohozaev15} requires that the solution $u$ is of class $C^2$. One can refer to the relevant literature in \cite{zengjia2,zengjia3} and the references therein.

Wagner \cite{Wanger} revisited the Pohozaev identity from the viewpoint of domain variations for a semilinear
elliptic equation, which is depends on the $C^2$ regularity of the solution. Xu \cite{Xu} presents a unified approach--domain deformation--to a class of identities which was first proved by Pohozaev, which also assumes that the solution is of class $C^2$. In \cite{Yan} by using the domain variations estimates, Yan and his collaborator derive the local Pohozaev identities for the weak solutions of p-Laplacian equation, without the $C^2$-regularity assumption for weak solutions. In this paper, we also obtain the Pohozaev identities for the weak solutions of Grushin type p-sub-Laplacian equation via domain variations, and we avoid the $C^2$-regularity assumption for weak solutions.

Now we provide an overview of the proof. Note that the definition \eqref{x420} of $g$. For any $g\in L^{p^*}(\Omega)$, where $p^*=\displaystyle\frac{p}{p-1}$, the problem \eqref{x6} has a unique solution $u\in W^{1,p}_{\lambda,0}(\Omega)$, which is proven in Lemma \ref{hg25}. Moreover, $u\in W^{1,p}_{\lambda,0}(\Omega)$ is the minimizer of the functional $I(u)$ defined by \eqref{x10}, which is proven in Lemma \ref{hg9}. Then, we will work on the functional $I(u)$. Let $\Phi_t(z):\Omega\rightarrow\Omega$ be such that $\Phi_0(z)=z$ and $\Phi_t(z):\partial\Omega\rightarrow\partial\Omega$. Then $u(\Phi_t(z))\in W^{1,p}_{\lambda,0}(\Omega)$ and $u(\Phi_0(z))=u(z)$. By the minimality of $u$, we know that
\begin{equation}\label{x11}
  \frac{d}{dt}I(u(\Phi_t(z)))\Big|_{t=0}=0.
\end{equation}
In the following, by choosing suitable $\Phi_t(z)$, we can compute each terms in \eqref{x11} to derive the local Pohozaev identities. In this process, we only need the $C^1$ regularity of weak solution $u$ and without any $C^2$ regularity assumption.

Our main results can be stated as follows. The following two lemmas establish the equivalence between the minimizer of  functional \eqref{x10} and the weak solution of equation \eqref{x6}.

\begin{Lem}\label{hg9}
For any $g\in L^{p^*}(\Omega)$, then there exists at least one function $u\in W^{1,p}_{\lambda,0}(\Omega)$ solving
\begin{equation}\label{hg10}
  I(u)=\inf_{W^{1,p}_{\lambda,0}(\Omega)}I(v).
\end{equation}
\end{Lem}

\begin{Lem}\label{hg25}
For any $g\in L^{p^*}(\Omega)$, there exists a unique weak solution for problem \eqref{x6}.
\end{Lem}

By working on functional $I(u)$, we derive the following Pohozaev identities.

\begin{Thm}\label{th1}\textbf{[Local Pohozaev identity of translating type]}\\
If $u\in W_{\lambda,0}^{1,p}(\Omega)\cap C^1(\Omega)$ is the solution of equation \eqref{x6}, then $u$ satisfies
\begin{equation}\label{Poho1}
\begin{aligned}
  \frac{1}{p}\int_{\partial D}|\nabla_\gamma u|^p\nu_x^idS
  -&\int_{\partial D}|\nabla_\gamma u|^{p-2}\frac{\partial u}{\partial x_i}\langle\nabla_\gamma u,\nu_\gamma\rangle dS
  -\int_D|\nabla_\gamma u|^{p-2}\gamma|x|^{2(\gamma-1)}x_i|\nabla_yu|^2dz
  \\=&\int_{\partial D}F(z,u)\nu_x^idS-\int_D\frac{\partial F(z,u)}{\partial x_i}dz,\;\;\;\;~\mbox~{i=1,\ldots,N},
\end{aligned}
\end{equation}
and
\begin{equation}\label{Poho2}
\begin{aligned}
  &\frac{1}{p}\int_{\partial D}|\nabla_\gamma u|^p\nu_y^jdS
  -\int_{\partial D}|\nabla_\gamma u|^{p-2}\frac{\partial u}{\partial y_j}\langle\nabla_\gamma u,\nu_\gamma\rangle dS
\\=&\int_{\partial D}F(z,u)\nu_y^jdS-\int_D\frac{\partial F(z,u)}{\partial y_j}(z,u)dz,\;\;\;\;~\mbox~{j=1,\ldots,l},
\end{aligned}
\end{equation}
where $D\subset\Omega$ is an any domain of $\Omega$, $\nu=(\nu_x,\nu_y)$ is the unit outward normal of the point of $\partial D$, $\nu_\gamma=(\nu_x,|x|^\gamma\nu_y)$ and $\displaystyle\frac{\partial F(z,u)}{\partial x_i}$ satisfies
\begin{equation}\label{sjm1}
  F_{x_i}(z,u)=\frac{\partial F(z,u)}{\partial x_i}+f(z,u)\frac{\partial u}{\partial x_i}.
\end{equation}
\end{Thm}

\begin{Thm}\label{th2}\textbf{[Local Pohozaev identity of scaling type]}\\
If $u\in W_{\lambda,0}^{1,p}(\Omega)\cap C^1(\Omega)$ is the solution of equation \eqref{x6}, then $u$ satisfies
\begin{equation}\label{Poho3}
\begin{aligned}
&\left(-\frac{N+l}{p}+1\right)\int_Df(z,u)udz+\left(-\frac{N+l}{p}+1\right)\int_{\partial D}|\nabla_\gamma u|^{p-2}u\langle\nabla_\gamma u,\nu_\gamma\rangle dS
 \\& +\frac{1}{p}\int_{\partial D}|\nabla_\gamma u|^p\langle z,\nu\rangle dS
-\int_{\partial D}|\nabla_\gamma u|^{p-2}\langle\nabla u,z\rangle\langle\nabla_\gamma u,\nu_\gamma\rangle dS
-\int_D|\nabla_\gamma u|^{p-2}\gamma|x|^{2\gamma}|\nabla_yu|^2dz
\\=&\int_{\partial D}F(z,u)\langle z,\nu\rangle dS
  -(N+l)\int_DF(z,u)dz
  -\int_D\langle \nabla_zF(z,u),z\rangle dz,
\end{aligned}
\end{equation}
where $D\subset\Omega$ is an any domain of $\Omega$, $\nu=(\nu_x,\nu_y)$ is the unit outward normal of the point of $\partial D$, $\nu_\gamma=(\nu_x,|x|^\gamma\nu_y)$, and $\nabla_zF(z,u)$ satisfies
\begin{equation}\label{sjm2}
  \nabla F(z,u)=\nabla_zF(z,u)+f(z,u)\nabla u.
\end{equation}
\end{Thm}

\begin{Thm}\label{th3}\textbf{[Global Pohozaev identity of scaling type]}\\
If $u\in W^{1,p}_{\lambda,0}(\Omega)\cap C^1(\overline{\Omega})$ is the solution of equation \eqref{x6}, then $u$ satisfies
\begin{equation}\label{d37}
\begin{aligned}
&\left(-\frac{N+l}{p}+1\right)\int_\Omega f(z,u)udz
+\left(\frac{1}{p}-1\right)\int_{\partial\Omega}|\nabla_\gamma u|^p\langle z,\nu\rangle dS
 \\-&\int_\Omega|\nabla_\gamma u|^{p-2}\gamma|x|^{2\gamma}|\nabla_yu|^2dz
=
  -(N+l)\int_\Omega F(z,u)dz
  -\int_\Omega\langle \nabla_zF(z,u),z\rangle dz,
\end{aligned}
\end{equation}
where $\Omega\subset\R^{N+l}$ is an any domain of $\R^{N+l}$, $\nu=(\nu_x,\nu_y)$ is the unit outward normal of the point of $\partial\Omega$, and $\nabla_z F(z,u)$ satisfies \eqref{sjm2}.

Moreover, if $\Omega=\R^{N+l}$, then $u$ satisfies
\begin{equation}\label{d38}
\begin{aligned}
&\left(-\frac{N+l}{p}+1\right)\int_{\R^{N+l}}f(z,u)udz
-\int_{\R^{N+l}}|\nabla_\gamma u|^{p-2}\gamma|x|^{2\gamma}|\nabla_yu|^2dz
\\=&
  -(N+l)\int_{\R^{N+l}}F(z,u)dz
  -\int_{\R^{N+l}}\langle \nabla_zF(z,u),z\rangle dz.
\end{aligned}
\end{equation}
\end{Thm}

The main contributions of this paper are summarized in the following three points. First, since there is no symmetric invariance for the $x$ variable in Grushin type p-sub-Laplacian \eqref{x5}, we establish two local Pohozaev identities of translating type by choosing $\Phi_t(z)=z+t\phi(z)e_i$ and $\Phi_t(z)=z+t\phi(z)e_{N+j}$ respectively for $1\leq i\leq N$ and $1\leq j\leq l$. Second, we obtain local Pohozaev identity of scaling type by choosing $\Phi_t(z)=(1+t\phi(z))z$ to overcome the challenge from the asymmetry of the Grushin operators. Third, as an application, a global Pohozaev identity of scaling type in $\R^{N+l}$ is also given, which can be used to construct Pohozaev manifolds.

The paper in the sequel is organized as follows. In Section 2, we prove Lemma \ref{hg9} and Lemma \ref{hg25} to obtain the equivalence between the minimizer of functional $I(u)$ and  the weak solution of equation \eqref{x6}. In Section 3, we prove Theorem \ref{th1} to obtain local Pohozaev identities generated from translations of \eqref{x6}. In Section 4, we prove Theorem \ref{th2} to obtain local Pohozaev identity generated from scaling of \eqref{x6}. In Section 5, we prove Theorem \ref{th3} to obtain a global Pohozaev identity of scaling type in $\R^{N+l}$.

\section{The equivalence between functional minimizer and weak solution}

In this section, in order to work on the functional $I(u)$, we give two lemms. In Lemma \ref{hg9}, we prove the existence of minimizer for the functional $I(u)$. In Lemma \ref{hg25}, we prove problem \eqref{x6} has a unique weak solution. It follows that seeking the weak solution of \eqref{x6} is equivalent to finding the minimizer of \eqref{x10}.

First, we present the Poincar\'{e} inequality and the Young's inequality with $\varepsilon$, which are crucial in the proofs of two lemmas.

\begin{Prop}\label{Pang}\textbf{(Degenerate Friedrichs-Poincar\'{e} inequality, \cite{chenhongge}).}\\
Suppose that $X=(X_1,X_2,\cdots,X_m)$ are the smooth vector fields defined on $U$, satisfying the H\"{o}rmander's condition. For any open bounded subset $W\subset\subset U$ and positive number $p\geq1$, there exists a positive constant $C>0$ such that
\begin{equation}\label{hg1}
  \int_W|u|^pdx\leq C\int_W|Xu|^pdx,\;\;\;\; \forall u\in \mathcal{W}^{1,p}_{X,0}(W).
\end{equation}
In which, $\mathcal{W}^{1,p}_{X,0}(W)$ is the closure of $C_0^\infty(W)$ in $\mathcal{W}^{1,p}_X(W)$ with respect to the norm $\|\cdot\|_{\mathcal{W}^{1,p}_X(W)}$.
For $p\geq1$,
\begin{equation}\label{hg2}
  \mathcal{W}^{1,p}_X(W)=\left\{u\in L^p(W)\;|\;X^Ju\in L^p(W), \;\;\forall J=(j_1,\cdots,j_m)~\mbox{with}~|J|\leq1.\right\},
\end{equation}
and
\begin{equation}\label{hg3}
  \|u\|^p_{\mathcal{W}^{1,p}_X(W)}=\sum_{|J|\leq1}\|X^Ju\|^p_{L^p(W)}.
\end{equation}
$\mathcal{W}^{1,p}_X(W)$ is a reflexive Banach space for $1<p<+\infty$.
\end{Prop}

\begin{Rem}\label{hg4}
Grushin operators form classical  examples in the family of all H\"ormander operators. They include the operator $\Delta_\gamma$ as an important model case. More generally, when we take
\begin{equation}\label{hg5}
  X=(\frac{\partial}{\partial x_1},\cdots,\frac{\partial}{\partial x_N},|x|^\gamma\frac{\partial}{\partial y_1},\cdots,|x|^\gamma\frac{\partial}{\partial y_l}),
\end{equation}
we can similarly obtain the weighted Poincar\'{e} inequality for the Grushin vector fields as follows (cf. \cite[Proposition 2.3]{chenhongge1}).
For any open bounded subset $\Omega\subset\R^{N+l}$ and positive number $p\geq1$, there exists a positive constant $\mu>0$ such that
\begin{equation}\label{hg6}
 \int_\Omega|u|^pdz\leq\mu\int_\Omega|\nabla_\gamma u|^pdz,\;\;\;\; \forall u\in W^{1,p}_{\lambda,0}(\Omega).
\end{equation}
Furthermore, $W^{1,p}_{\lambda,0}(\Omega)$ and $W^{1,p}_\lambda(\Omega)$ defined in \eqref{x8} are reflexive Banach spaces.
\end{Rem}

\begin{Prop}\label{hg7}\textbf{(Young's inequality with $\varepsilon$, \cite{Evans}).}\\
Let $a>0$, $b>0$, $\varepsilon>0$, $p>1$, $p^*>1$, and $\displaystyle\frac{1}{p}+\frac{1}{p^*}=1$, then
\begin{equation}\label{hg8}
  ab\leq\frac{\varepsilon}{p}a^p+\frac{1}{p^*\varepsilon^{\frac{p^*}{p}}}b^{p^*}.
\end{equation}
\end{Prop}

Now, we present the proofs for Lemma \ref{hg9} and Lemma \ref{hg25}.

\begin{proof}[\textbf{Proof of Lemma \ref{hg9}:}]

\vskip 0.2cm

\noindent\textbf{Step 1: To prove the existence of minimizing sequence} \{$u_k$\} \textbf{for functional $I(u)$ in $W^{1,p}_{\lambda,0}(\Omega)$}.

\vskip 0.2cm

By Poincar\'{e} inequality \eqref{hg6} and Young's inequality \eqref{hg8}, we know that
\begin{equation}\label{hg11}
\begin{aligned}
 I(u)&=\frac{1}{p}\int_\Omega|\nabla_\gamma u|^pdz-\int_\Omega g(z)udz
 \\&\geq\frac{1}{p\mu}\int_\Omega|u|^pdz-\int_\Omega\left(\frac{\varepsilon}{p}|u|^p
 +\frac{1}{p^*\varepsilon^{\frac{p^*}{p}}}|g|^{p^*}\right)dz
 \\&=\frac{1}{p}\left(\frac{1}{\mu}-\varepsilon\right)\int_\Omega|u|^pdz
 -\frac{1}{p^*\varepsilon^{\frac{p^*}{p}}}\int_\Omega|g|^{p^*}dz,
\end{aligned}
\end{equation}
where $\mu>0$ is a constant in \eqref{hg6}, and $\varepsilon>0$ is a constant that can be taken arbitrarily. We take $\varepsilon\leq\frac{1}{\mu}$, then \eqref{hg11} turns into
\begin{equation}\label{hg12}
  I(u)\geq-\frac{1}{p^*\varepsilon^{\frac{p^*}{p}}}\int_\Omega|g|^{p^*}dz,
\end{equation}
which implies that $I(u)$ is bounded from below in $W^{1,p}_{\lambda,0}(\Omega)$ and thus it has infimum in $W^{1,p}_{\lambda,0}(\Omega)$.
Therefore by the definition of infimum, there exists $u_k\in W^{1,p}_{\lambda,0}(\Omega)$, such that
\begin{equation}\label{hg13}
  \lim_{k\rightarrow\infty}I(u_k)=\inf_{W^{1,p}_{\lambda,0}(\Omega)}I(v).
\end{equation}

\vskip 0.2cm

\noindent\textbf{Step 2: To prove there exist subsequence} \{$u_{k_i}$\}$\subset$\{$u_k$\} \textbf{and} $u\in W^{1,p}_{\lambda,0}(\Omega)$, \textbf{such that}
\begin{equation}\label{hg14}
  u_{k_i}\rightharpoonup u,\;\;~\mbox{weakly in}~W^{1,p}_{\lambda,0}(\Omega).
\end{equation}

\vskip 0.2cm

By \eqref{hg13}, the limit $\displaystyle\lim_{k\rightarrow\infty}I(u_k)$ exists, then there is a constant $M>0$, such that
\begin{equation}\label{hg15}
  |I(u_k)|\leq M.
\end{equation}
According to inequalities \eqref{hg6} and \eqref{hg8}, it holds
\begin{equation}\label{hg16}
\begin{aligned}
  \int_\Omega|\nabla_\gamma u_k|^pdz=&p\int_\Omega g(z)u_kdz+pI(u_k)
  \\\overset{\eqref{hg8}}\leq&\varepsilon\int_\Omega|u_k|^pdz
  +\frac{p}{p^*\varepsilon^{\frac{p^*}{p}}}\int_\Omega|g|^{p^*}dz+pM
   \\\overset{\eqref{hg6}}\leq&\varepsilon\mu\int_\Omega|\nabla_\gamma u_k|^pdz
  +\frac{p}{p^*\varepsilon^{\frac{p^*}{p}}}\int_\Omega|g|^{p^*}dz+pM.
\end{aligned}
\end{equation}
Take $\varepsilon=\displaystyle\frac{1}{2\mu}$ in \eqref{hg16}, we obtain
\begin{equation}\label{hg17}
  \int_\Omega|\nabla_\gamma u_k|^pdz\leq\frac{2 p}{p^*}(2\mu)^{\frac{p^*}{p}}\int_\Omega|g|^{p^*}dz+2pM.
\end{equation}
From \eqref{hg17} and by inequality \eqref{hg6} again, we obtain
\begin{equation}\label{hg18}
  \int_\Omega|u_k|^pdz\leq\frac{p}{p^*}(2\mu)^{p^*}\int_\Omega|g|^{p^*}dz+2\mu pM.
\end{equation}
By \eqref{hg17} and \eqref{hg18}, we know that \{$u_k$\} is bounded in $W^{1,p}_{\lambda,0}(\Omega)$. Since $W^{1,p}_{\lambda,0}(\Omega)$ is a reflexive Banach space, then there exist subsequence \{$u_{k_i}$\}$\subset$\{$u_k$\} and $u\in W^{1,p}_{\lambda,0}(\Omega)$, satisfying \eqref{hg14}.

\vskip 0.2cm

\noindent\textbf{Step 3: To prove $I(u)$ is weakly lower semicontinuous on $W^{1,p}_{\lambda,0}(\Omega)$}.

\vskip 0.2cm

For any \{$u_k$\}$\subset$$W^{1,p}_{\lambda,0}(\Omega)$, satisfying
\begin{equation}\label{hg19}
    u_k\rightharpoonup u,\;\;~\mbox{weakly in}~W^{1,p}_{\lambda,0}(\Omega),
\end{equation}
we want to prove
\begin{equation}\label{hg20}
  I(u)\leq\liminf_{k\rightarrow\infty}I(u_k).
\end{equation}
By the weakly lower semicontinuity of the $L^p$ norm, we have
\begin{equation}\label{hg21}
  \int_\Omega|\nabla_\gamma u|^pdz\leq\liminf_{k\rightarrow\infty}\int_\Omega|\nabla_\gamma u_k|^pdz.
\end{equation}
Since \eqref{hg19}, it holds that
\begin{equation}\label{hg22}
  \lim_{k\rightarrow\infty}\int_\Omega g(z)u_kdz=\int_\Omega g(z)udz.
\end{equation}
According to \eqref{hg21} and \eqref{hg22}, we conclude that
\begin{equation}\label{hg23}
\begin{aligned}
  \liminf_{k\rightarrow\infty}I(u_k)&=\frac{1}{p}\liminf_{k\rightarrow\infty}\int_\Omega|\nabla_\gamma u_k|^pdz
  -\lim_{k\rightarrow\infty}\int_\Omega g(z)u_kdz
  \\&\geq\frac{1}{p}\int_\Omega|\nabla_\gamma u|^pdz-\int_\Omega g(z)udz
  \\&=I(u),
\end{aligned}
\end{equation}
which means that $I(u)$ is weakly lower semicontinuous on $W^{1,p}_{\lambda,0}(\Omega)$.

\vskip 0.2cm

\noindent\textbf{Step 4: Complete the proof}.

\vskip 0.2cm

From \textbf{Step 1}--\textbf{Step 3}, we derive that
\begin{equation}\label{hg24}
  \inf_{W^{1,p}_{\lambda,0}(\Omega)}I(v)\leq I(u)\overset{\eqref{hg20}}\leq\liminf_{k\rightarrow\infty}I(u_{k_i})
  \overset{\eqref{hg13}}=\lim_{k\rightarrow\infty}I(u_{k_i})\overset{\eqref{hg13}}=\inf_{W^{1,p}_{\lambda,0}(\Omega)}I(v).
\end{equation}
Then we have
\begin{equation*}
  I(u)=\inf_{W^{1,p}_{\lambda,0}(\Omega)}I(v).
\end{equation*}

\end{proof}

\begin{proof}[\textbf{Proof of Lemma \ref{hg25}:}]
From Lemma \ref{hg9}, we know that the existence of the minimizer of functional $I(u)$. According to basic ideas of the calculus of variations, if $u$ is the minimizer of  $I(u)$, then $u$ is the solution of \eqref{x6}. Therefore, the existence of weak solution is obtained. Now, we prove the uniqueness of weak solution.

Let $u_1,u_2\in W^{1,p}_{\lambda,0}(\Omega)$ be two weak solutions of equation \eqref{x6}. By the definition of \eqref{x7}, we have
\begin{equation}\label{hg26}
   \int_\Omega|\nabla_\gamma u_i|^{p-2}\nabla_\gamma u_i\cdot\nabla_\gamma\varphi dz=\int_\Omega g(z)\varphi dz,\;\;\forall\varphi\in W_{\lambda,0}^{1,p}(\Omega),\;\;i=1,2.
\end{equation}
From \eqref{hg26}, we obtain
\begin{equation}\label{reject1}
 \int_\Omega\left(|\nabla_\gamma u_1|^{p-2}\nabla_\gamma u_1-|\nabla_\gamma u_2|^{p-2}\nabla_\gamma u_2\right)\cdot\nabla_\gamma\varphi dz=0,\;\;\forall\varphi\in W_{\lambda,0}^{1,p}(\Omega).
\end{equation}
Especially take $\varphi$ satisfying
\begin{equation}\label{reject2}
  \nabla_\gamma\varphi=|\nabla_\gamma u_1|^{p-2}\nabla_\gamma u_1+|\nabla_\gamma u_2|^{p-2}\nabla_\gamma u_2,
\end{equation}
then \eqref{reject1} transforms into
\begin{equation}\label{reject3}
   \int_\Omega|\nabla_\gamma u_1|^{2p-2}-|\nabla_\gamma u_2|^{2p-2}dz=0.
\end{equation}
Since $p>1$, $2p-2>0$, then from \eqref{reject3} we have
\begin{equation}\label{reject3}
   \int_\Omega|\nabla_\gamma u_1|^p-|\nabla_\gamma u_2|^pdz=0.
\end{equation}

On the other hand, take $\varphi=u_1$ and $\varphi=u_2$ respectively in \eqref{hg26}, it holds that
\begin{equation}\label{reject4}
  \int_\Omega|\nabla_\gamma u_i|^pdz=\int_\Omega g(z)u_i dz,\;\;i=1,2.
\end{equation}
From \eqref{reject3} and \eqref{reject4}, we have
\begin{equation}\label{reject5}
  \int_\Omega g(z)\left(u_1-u_2\right)dz=0.
\end{equation}
For any $g\in L^{p^*}(\Omega)$, it follows from \eqref{reject5} that
\begin{equation}\label{reject5}
  u_1=u_2,\;\;\;~\mbox{a.e. in}~\Omega.
\end{equation}

\end{proof}

\section{Local Pohozaev identity of translating type}

In this section, we prove two Pohozaev identities generated from translations of \eqref{x6} in Theorem \ref{th1}.

\begin{proof}[\textbf{Proof of Theorem \ref{th1}:}]
Suppose that $D\subset\subset\Omega$, and for any $\delta>0$ small enough define
\begin{equation}\label{x14}
  D_\delta:=\{z\in\Omega\;|\;d(z,D)<\delta\}\subset\Omega.
\end{equation}
Now, we take
\begin{equation}\label{x12}
  \Phi_t(z)=z+t\phi(z)e_j,
\end{equation}
where $t\in\R$, $z=(x,y)=(x_1,\cdots,x_N,y_1,\cdots,y_l)\in\Omega$, $e_j =(\underbrace{0, \ldots, 1}_{j}, \ldots, 0)$, $\phi\in W^{1,p}_{\lambda,0}(\Omega)$ satisfies
\begin{equation}\label{x13}
\phi(z)=\left\{\begin{array}{ll}
1, &~\mbox{in}~D ,\\[2mm]
0, &~\mbox{in}~\Omega\setminus D_\delta,
\end{array}\right.
\end{equation}
and
\begin{equation}\label{x15}
  0\leq\phi(z)\leq1,~\mbox{in}~D_\delta\setminus D.
\end{equation}
Here, $\nu=(\nu_x,\nu_y)$ is the outward unit normal of $\partial D$ at $\overline{z}=(\overline{x},\overline{y})\in\partial D$ with
\begin{equation}\label{x16}
  d(z,D)=\left(\frac{1}{(1+\gamma)^2}|x-\overline{x}|^{2+2\gamma}+|y-\overline{y}|^2\right)^{\frac{1}{2+2\gamma}}.
\end{equation}
Note that
\begin{equation}\label{x38}
  \nabla\phi(z)=-\nu,\;\;\nabla_\gamma\phi(z)=-\nu_\gamma=-(\nu_x,|x|^\gamma\nu_y) ~\mbox{on}~ \partial D,\;\;~\mbox{as}~\delta\rightarrow0.
\end{equation}
Clearly, $\Phi_t(z):\Omega\rightarrow\Omega$ satisfies $\Phi_0(z)=z$ and $\Phi_t(z):\partial\Omega\rightarrow\partial\Omega$.

By the definition \eqref{x10} of functional $I$, we have
\begin{equation}\label{x180}
  I(u(\Phi_t(z)))=\frac{1}{p}\int_\Omega|\nabla_\gamma u(\Phi_t(z))|^pdz-\int_\Omega g(z)u(\Phi(z))dz.
\end{equation}
It follows from \eqref{x11} that
\begin{equation}\label{x17}
  \frac{d}{dt}I(u(\Phi_t(z)))\Big|_{t=0}=\frac{1}{p}\frac{d}{dt}\int_\Omega|\nabla_\gamma u(\Phi_t(z))|^pdz\Big|_{t=0}
  -\frac{d}{dt}\int_\Omega g(z)u(\Phi_t(z))dz\Big|_{t=0}=0.
\end{equation}
Without loss of generality, in what follows, we simply take $j=1$ and $j=N+1$ respectively.

\vskip 0.2cm

\textbf{We first consider the case of $j=1$}.

\vskip 0.2cm

\noindent\textbf{Step 1: Calculate} $\displaystyle\frac{d}{dt}~\mbox{det}~(Jac\Phi^{-1}_t(z))\Big|_{t=0}$.

\vskip 0.2cm

\begin{equation}\label{x18}
\begin{aligned}
  \Phi_t(z)&=z+t\phi(z)e_1
  \\&=(x_1+t\phi(z),x_2,\cdots,x_N,y_1,\cdots, y_l).
\end{aligned}
\end{equation}
Let
\begin{equation}\label{x20}
  \widetilde{z}=\Phi_t(z),
\end{equation}
i.e.
\begin{equation}\label{x19}
  \widetilde{x_1}=x_1+t\phi(z),\;\widetilde{x_2}=x_2,\;\cdots,\;\widetilde{x_N}=x_N,\;\widetilde{y_1}=y_1,\;\cdots,\; \widetilde{y_l}=y_l.
\end{equation}
By directly calculation, the Jacobian matrix of $\Phi_t(z)$ is as follows
\begin{equation}\label{x28}
Jac\Phi_t(z)=
\begin{bmatrix}
\frac{\partial\widetilde{x_1}}{\partial x_1} & \frac{\partial\widetilde{x_1}}{\partial x_2} & \cdots & \frac{\partial\widetilde{x_1}}{\partial y_l}\\
\frac{\partial\widetilde{x_2}}{\partial x_1} & \frac{\partial\widetilde{x_2}}{\partial x_2} & \cdots & \frac{\partial\widetilde{x_2}}{\partial y_l}\\
\vdots & \vdots & \ddots & \vdots\\
\frac{\partial\widetilde{y_l}}{\partial x_1} & \frac{\partial\widetilde{y_l}}{\partial x_2} & \cdots & \frac{\partial\widetilde{y_l}}{\partial y_l}\\
\end{bmatrix}
=
\begin{bmatrix}
1+t\frac{\partial\phi(z)}{\partial x_1} & t\frac{\partial\phi(z)}{\partial x_2} & \cdots & t\frac{\partial\phi(z)}{\partial y_l}\\
0 & 1 & \cdots & 0\\
\vdots & \vdots & \ddots & \vdots\\
0 & 0 & \cdots & 1\\
\end{bmatrix}.
\end{equation}
From \eqref{x28}, we derive that
\begin{equation}\label{x29}
  ~\mbox{det}~Jac\Phi_t(z)=1+t\frac{\partial\phi(z)}{\partial x_1}.
\end{equation}
Since
\begin{equation}\label{x30}
 ~\mbox{det}~Jac\Phi_t(z)\cdot~\mbox{det}~Jac\Phi_t^{-1}(z)=1,
\end{equation}
then
\begin{equation}\label{x31}
  ~\mbox{det}~Jac\Phi_t^{-1}(z)=\frac{1}{1+t\frac{\partial\phi(z)}{\partial x_1}}.
\end{equation}
Thus we obtain
\begin{equation}\label{x32}
  \frac{d}{dt}~\mbox{det}~Jac\Phi_t^{-1}(z)=-\frac{\frac{\partial\phi(z)}{\partial x_1}}{(1+t\frac{\partial\phi(z)}{\partial x_1})^2}.
\end{equation}
Accordingly, we have
\begin{equation}\label{x33}
  \frac{d}{dt}~\mbox{det}~(Jac\Phi^{-1}_t(z))\Big|_{t=0}=-\frac{\partial\phi(z)}{\partial x_1}.
\end{equation}

\vskip 0.2cm

\noindent\textbf{Step 2: Calculate $\displaystyle\frac{1}{p}\frac{d}{dt}\int_\Omega|\nabla_\gamma u(\Phi_t(z))|^pdz\Big|_{t=0}$.}

\vskip 0.2cm

Note that
\begin{equation}\label{x35}
  \frac{d}{dt}\int_\Omega|\nabla_\gamma u(\Phi_t(z))|^pdz\Big|_{t=0}=\int_\Omega\frac{d}{dt}|\nabla_\gamma u(\Phi_t(z))|^p\Big|_{t=0}dz.
\end{equation}
According to \eqref{x18}-\eqref{x19}, by calculation we have
\begin{equation}\label{x21}
  |\nabla_xu(\widetilde{z})|^2=\sum_{i=1}^N\left(\frac{\partial u(\widetilde{z})}{\partial\widetilde{x_1}}
  t\frac{\partial\phi(z)}{\partial x_i}+\frac{\partial u(\widetilde{z})}{\partial\widetilde{x_i}}\right)^2,
\end{equation}
\begin{equation}\label{x22}
  |\nabla_yu(\widetilde{z})|^2=\sum_{j=1}^l\left(\frac{\partial u(\widetilde{z})}{\partial\widetilde{x_1}}
  t\frac{\partial\phi(z)}{\partial y_j}+\frac{\partial u(\widetilde{z})}{\partial\widetilde{y_j}}\right)^2,
\end{equation}
and
\begin{equation}\label{x23}
  |x|^{2\gamma}=\left(|\widetilde{x}|^2+t^2\phi^2(z)-2\widetilde{x_1}t\phi(z)\right)^\gamma.
\end{equation}
Thus, it holds
\begin{equation}\label{x24}
\begin{aligned}
  |\nabla_\gamma u(\widetilde{z})|^2=&\sum_{i=1}^N\left(\frac{\partial u(\widetilde{z})}{\partial\widetilde{x_1}}
  t\frac{\partial\phi(z)}{\partial x_i}+\frac{\partial u(\widetilde{z})}{\partial\widetilde{x_i}}\right)^2
  \\+&\left(|\widetilde{x}|^2+t^2\phi^2(z)-2\widetilde{x_1}t\phi(z)\right)^\gamma
  \sum_{j=1}^l\left(\frac{\partial u(\widetilde{z})}{\partial\widetilde{x_1}}
  t\frac{\partial\phi(z)}{\partial y_j}+\frac{\partial u(\widetilde{z})}{\partial\widetilde{y_j}}\right)^2,
\end{aligned}
\end{equation}
and
\begin{equation}\label{x25}
\begin{aligned}
  |\nabla_\gamma u(\widetilde{z})|^p&=\left(\sum_{i=1}^N\left(\frac{\partial u(\widetilde{z})}{\partial\widetilde{x_1}}
  t\frac{\partial\phi(z)}{\partial x_i}+\frac{\partial u(\widetilde{z})}{\partial\widetilde{x_i}}\right)^2\right.
  \\&\left.+\left(|\widetilde{x}|^2+t^2\phi^2(z)-2\widetilde{x_1}t\phi(z)\right)^\gamma
  \sum_{j=1}^l\left(\frac{\partial u(\widetilde{z})}{\partial\widetilde{x_1}}
  t\frac{\partial\phi(z)}{\partial y_j}+\frac{\partial u(\widetilde{z})}{\partial\widetilde{y_j}}\right)^2\right)^{\frac{p}{2}}.
\end{aligned}
\end{equation}
Accordingly, we have
\begin{equation}\label{x26}
\begin{aligned}
  &\int_\Omega|\nabla_\gamma u(\Phi_t(z))|^pdz
  \\=&\int_\Omega\left(\sum_{i=1}^N\left(\frac{\partial u(\widetilde{z})}{\partial\widetilde{x_1}}
  t\frac{\partial\phi(z)}{\partial x_i}+\frac{\partial u(\widetilde{z})}{\partial\widetilde{x_i}}\right)^2\right.
  \\&\left.+\left(|\widetilde{x}|^2+t^2\phi^2(z)-2\widetilde{x_1}t\phi(z)\right)^\gamma
  \sum_{j=1}^l\left(\frac{\partial u(\widetilde{z})}{\partial\widetilde{x_1}}
  t\frac{\partial\phi(z)}{\partial y_j}+\frac{\partial u(\widetilde{z})}{\partial\widetilde{y_j}}\right)^2\right)^{\frac{p}{2}}dz
  \\=&\int_\Omega\left(\sum_{i=1}^N\left(\frac{\partial u(\widetilde{z})}{\partial\widetilde{x_1}}
  t\frac{\partial\phi(\Phi_t^{-1}(\widetilde{z}))}{\partial (\Phi_t^{-1}(\widetilde{z}))_i}+\frac{\partial u(\widetilde{z})}{\partial\widetilde{x_i}}\right)^2\right.
  \\&\left.+\left(|\widetilde{x}|^2+t^2\phi^2(\Phi_t^{-1}(\widetilde{z}))-2\widetilde{x_1}
  t\phi(\Phi_t^{-1}(\widetilde{z}))\right)^\gamma
  \sum_{j=1}^l\left(\frac{\partial u(\widetilde{z})}{\partial\widetilde{x_1}}
  t\frac{\partial\phi(\Phi_t^{-1}(\widetilde{z}))}{\partial (\Phi_t^{-1}(\widetilde{z}))_{N+j}}+\frac{\partial u(\widetilde{z})}{\partial\widetilde{y_j}}\right)^2\right)^{\frac{p}{2}}
  \\&~\mbox{det}~(Jac\Phi^{-1}_t(\widetilde{z}))d\widetilde{z}
  \\=&\int_\Omega\left(\sum_{i=1}^N\left(\frac{\partial u(z)}{\partial x_1}
  t\frac{\partial\phi(\Phi_t^{-1}(z))}{\partial (\Phi_t^{-1}(z))_i}+\frac{\partial u(z)}{\partial x_i}\right)^2\right.
  \\&\left.+\left(|x|^2+t^2\phi^2(\Phi_t^{-1}(z))-2x_1
  t\phi(\Phi_t^{-1}(z))\right)^\gamma
  \sum_{j=1}^l\left(\frac{\partial u(z)}{\partial x_1}
  t\frac{\partial\phi(\Phi_t^{-1}(z))}{\partial (\Phi_t^{-1}(z))_{N+j}}+\frac{\partial u(z)}{\partial y_j}\right)^2\right)^{\frac{p}{2}}
  \\&~\mbox{det}~(Jac\Phi^{-1}_t(z))dz.
\end{aligned}
\end{equation}
Then by \eqref{x35}, we have

\begin{equation*}
\begin{aligned}
 &\frac{d}{dt}\int_\Omega|\nabla_\gamma u(\Phi_t(z))|^pdz
\\=&\int_\Omega\frac{d}{dt}~\mbox{det}~(Jac\Phi^{-1}_t(z))
\cdot\left(\sum_{i=1}^N\left(\frac{\partial u(z)}{\partial x_1}
  t\frac{\partial\phi(\Phi_t^{-1}(z))}{\partial (\Phi_t^{-1}(z))_i}+\frac{\partial u(z)}{\partial x_i}\right)^2\right.
  \\&\left.+\left(|x|^2+t^2\phi^2(\Phi_t^{-1}(z))-2x_1
  t\phi(\Phi_t^{-1}(z))\right)^\gamma
  \sum_{j=1}^l\left(\frac{\partial u(z)}{\partial x_1}
  t\frac{\partial\phi(\Phi_t^{-1}(z))}{\partial (\Phi_t^{-1}(z))_{N+j}}+\frac{\partial u(z)}{\partial y_j}\right)^2\right)^{\frac{p}{2}}dz
\end{aligned}
\end{equation*}

\begin{equation}\label{x27}
\begin{aligned}
+&\int_\Omega~\mbox{det}~(Jac\Phi^{-1}_t(z))
\cdot\frac{p}{2}\left(\sum_{i=1}^N\left(\frac{\partial u(z)}{\partial x_1}
  t\frac{\partial\phi(\Phi_t^{-1}(z))}{\partial (\Phi_t^{-1}(z))_i}+\frac{\partial u(z)}{\partial x_i}\right)^2\right.
  \\&\left.+\left(|x|^2+t^2\phi^2(\Phi_t^{-1}(z))-2x_1
  t\phi(\Phi_t^{-1}(z))\right)^\gamma
  \sum_{j=1}^l\left(\frac{\partial u(z)}{\partial x_1}
  t\frac{\partial\phi(\Phi_t^{-1}(z))}{\partial (\Phi_t^{-1}(z))_{N+j}}+\frac{\partial u(z)}{\partial y_j}\right)^2\right)^{\frac{p}{2}-1}
\\\cdot&\left(\sum_{i=1}^N2\left(\frac{\partial u(z)}{\partial x_1}
  t\frac{\partial\phi(\Phi_t^{-1}(z))}{\partial (\Phi_t^{-1}(z))_i}+\frac{\partial u(z)}{\partial x_i}\right)\frac{\partial u(z)}{\partial x_1}
  \frac{\partial\phi(\Phi_t^{-1}(z))}{\partial (\Phi_t^{-1}(z))_i}\right.
  \\&+\left.\gamma\left(|x|^2+t^2\phi^2(\Phi_t^{-1}(z))-2x_1
  t\phi(\Phi_t^{-1}(z))\right)^{\gamma-1}(2t\phi^2(\Phi_t^{-1}(z)-2x_1\phi(\Phi_t^{-1}(z)))\right.
  \\&\cdot\left.\sum_{j=1}^l\left(\frac{\partial u(z)}{\partial x_1}
  t\frac{\partial\phi(\Phi_t^{-1}(z))}{\partial (\Phi_t^{-1}(z))_{N+j}}+\frac{\partial u(z)}{\partial y_j}\right)^2\right.
  +\left(|x|^2+t^2\phi^2(\Phi_t^{-1}(z))-2x_1
  t\phi(\Phi_t^{-1}(z))\right)^\gamma
  \\&\cdot\left.\sum_{j=1}^l2\left(\frac{\partial u(z)}{\partial x_1}
  t\frac{\partial\phi(\Phi_t^{-1}(z))}{\partial (\Phi_t^{-1}(z))_{N+j}}+\frac{\partial u(z)}{\partial y_j}\right)\frac{\partial u(z)}{\partial x_1}
  \frac{\partial\phi(\Phi_t^{-1}(z))}{\partial (\Phi_t^{-1}(z))_{N+j}}\right)dz.
\end{aligned}
\end{equation}

Since $\Phi_0(z)=z$, then $\Phi_0^{-1}(z)=z$ and $~\mbox{det}~Jac\Phi_0^{-1}(z)=1$. And by \textbf{Step 1}, we get

\begin{equation}\label{x34}
\begin{aligned}
 &\frac{d}{dt}\int_\Omega|\nabla_\gamma u(\Phi_t(z))|^pdz\Big|_{t=0}
\\=&\int_\Omega-\frac{\partial\phi(z)}{\partial x_1}\left(\sum_{i=1}^N\left(\frac{\partial u(z)}{\partial x_i}\right)^2+|x|^{2\gamma}
  \sum_{j=1}^l\left(\frac{\partial u(z)}{\partial y_j}\right)^2\right)^{\frac{p}{2}}
\\&+\frac{p}{2}\left(\sum_{i=1}^N\left(\frac{\partial u(z)}{\partial x_i}\right)^2+|x|^{2\gamma}
  \sum_{j=1}^l\left(\frac{\partial u(z)}{\partial y_j}\right)^2\right)^{\frac{p}{2}-1}
\\\cdot&\left(\sum_{i=1}^N2\frac{\partial u(z)}{\partial x_i}\frac{\partial u(z)}{\partial x_1}
  \frac{\partial\phi(z)}{\partial x_i}
  +\gamma|x|^{2(\gamma-1)}(-2x_1\phi(z))\sum_{j=1}^l\left(\frac{\partial u(z)}{\partial y_j}\right)^2
  +|x|^{2\gamma}
  \sum_{j=1}^l2\frac{\partial u(z)}{\partial y_j}\frac{\partial u(z)}{\partial x_1}
  \frac{\partial\phi(z)}{\partial y_j}\right)dz
\\=&\int_\Omega-\frac{\partial\phi(z)}{\partial x_1}|\nabla_\gamma u|^p+p|\nabla_\gamma u|^{p-2}\frac{\partial u}{\partial x_1}\langle\nabla_\gamma u,\nabla_\gamma \phi\rangle-p|\nabla_\gamma u|^{p-2}\gamma|x|^{2(\gamma-1)}x_1\phi(z)|\nabla_yu|^2dz.
\end{aligned}
\end{equation}
In which, by \eqref{x13} and \eqref{x38} we obtain
\begin{equation}\label{x37}
  \int_\Omega-\frac{\partial\phi(z)}{\partial x_1}|\nabla_\gamma u|^pdz
 =\int_{D_\delta\setminus D}-\frac{\partial\phi(z)}{\partial x_1}|\nabla_\gamma u|^pdz
 \overset{\delta\rightarrow0}\longrightarrow\int_{\partial D}|\nabla_\gamma u|^p\nu_x^1dS,
\end{equation}

\begin{equation}\label{x39}
\begin{aligned}
  &\int_\Omega|\nabla_\gamma u|^{p-2}\frac{\partial u}{\partial x_1}\langle\nabla_\gamma u,\nabla_\gamma \phi\rangle dz
\\=&\int_{D_\delta\setminus D}|\nabla_\gamma u|^{p-2}\frac{\partial u}{\partial x_1}\langle\nabla_\gamma u,\nabla_\gamma \phi\rangle dz
\\\overset{\delta\rightarrow0}\longrightarrow&-\int_{\partial D}|\nabla_\gamma u|^{p-2}\frac{\partial u}{\partial x_1}\langle\nabla_\gamma u,\nu_\gamma\rangle dS,
\end{aligned}
\end{equation}
and
\begin{equation}\label{x51}
\begin{aligned}
  &\int_\Omega|\nabla_\gamma u|^{p-2}\gamma|x|^{2(\gamma-1)}x_1\phi(z)|\nabla_yu|^2dz
  \\=&\int_{D_\delta}|\nabla_\gamma u|^{p-2}\gamma|x|^{2(\gamma-1)}x_1\phi(z)|\nabla_yu|^2dz
  \\\overset{\delta\rightarrow0}\longrightarrow&\int_D|\nabla_\gamma u|^{p-2}\gamma|x|^{2(\gamma-1)}x_1|\nabla_yu|^2dz.
\end{aligned}
\end{equation}

Therefore, by \eqref{x34}-\eqref{x51} we have
\begin{equation}\label{x40}
\begin{aligned}
  &\frac{1}{p}\frac{d}{dt}\int_\Omega|\nabla_\gamma u(\Phi_t(z))|^pdz\Big|_{t=0}
  \\\overset{\delta\rightarrow0}\longrightarrow&\frac{1}{p}\int_{\partial D}|\nabla_\gamma u|^p\nu_x^1dS
  -\int_{\partial D}|\nabla_\gamma u|^{p-2}\frac{\partial u}{\partial x_1}\langle\nabla_\gamma u,\nu_\gamma\rangle dS
  -\int_D|\nabla_\gamma u|^{p-2}\gamma|x|^{2(\gamma-1)}x_1|\nabla_yu|^2dz.
\end{aligned}
\end{equation}

\vskip 0.2cm

\noindent\textbf{Step 3: Calculate $\displaystyle\frac{d}{dt}\int_\Omega g(z)u(\Phi_t(z))dz\Big|_{t=0}$.}

\vskip 0.2cm

Note that
\begin{equation}\label{x41}
  \frac{d}{dt}\int_\Omega g(z)u(\Phi_t(z))dz\Big|_{t=0}
  =\int_\Omega g(z)\frac{d}{dt}u(\Phi_t(z))\Big|_{t=0}dz.
\end{equation}
Since
\begin{equation}\label{x42}
  u(\Phi_t(z))=u(z+t\phi(z)e_1),
\end{equation}
and by \eqref{x19}, it is easy to find that
\begin{equation}\label{x43}
  \frac{d}{dt}u(\Phi_t(z))=\frac{\partial u(\Phi_t(z))}{\partial\widetilde{x_1}}\phi(z),
\end{equation}
and thus
\begin{equation}\label{x44}
  \frac{d}{dt}u(\Phi_t(z))\Big|_{t=0}=\frac{\partial u(z)}{\partial x_1}\phi(z).
\end{equation}
Then
\begin{equation}\label{x45}
\begin{aligned}
  &\int_\Omega g(z)\frac{d}{dt}u(\Phi_t(z))\Big|_{t=0}dz
\\=&\int_\Omega g(z)\frac{\partial u(z)}{\partial x_1}\phi(z)dz
\\=&\int_{D_\delta} g(z)\frac{\partial u(z)}{\partial x_1}\phi(z)dz
\\\overset{\delta\rightarrow0}\longrightarrow&\int_D g(z)\frac{\partial u(z)}{\partial x_1}dz.
\end{aligned}
\end{equation}
On the other hand, by the definition \eqref{x420} of $g(z)$, we have
\begin{equation}\label{x46}
\begin{aligned}
  &\int_D g(z)\frac{\partial u(z)}{\partial x_1}dz
  \\=&\int_D f(z,u)\frac{\partial u(z)}{\partial x_1}dz
  \\=&\int_D \frac{\partial F(z,u)}{\partial x_1}-F_{x_1}(z,u)dz
  \\=&\int_{\partial D}F(z,u)\nu_x^1dS-\int_D\frac{\partial F(z,u)}{\partial x_1}dz,
\end{aligned}
\end{equation}
where $\displaystyle\frac{\partial F(z,u)}{\partial x_1}$ satisfies
\begin{equation}\label{sjm11}
  F_{x_1}(z,u)=\frac{\partial F(z,u)}{\partial x_1}+f(z,u)\frac{\partial u}{\partial x_1}.
\end{equation}
Thus from \eqref{x41}, \eqref{x45} and \eqref{x46}, we have
\begin{equation}\label{x48}
  \frac{d}{dt}\int_\Omega g(z)u(\Phi_t(z))dz\Big|_{t=0}\overset{\delta\rightarrow0}\longrightarrow\int_{\partial D}F(z,u)\nu_x^1dS-\int_D\frac{\partial F(z,u)}{\partial x_1}dz.
\end{equation}

\vskip 0.2cm

\noindent\textbf{Step 4: Complete the proof.}

\vskip 0.2cm

According to \eqref{x17}, \eqref{x40} and \eqref{x48}, we have
\begin{equation}\label{x49}
\begin{aligned}
  \frac{1}{p}\int_{\partial D}|\nabla_\gamma u|^p\nu_x^1dS
  -&\int_{\partial D}|\nabla_\gamma u|^{p-2}\frac{\partial u}{\partial x_1}\langle\nabla_\gamma u,\nu_\gamma\rangle dS
  -\int_D|\nabla_\gamma u|^{p-2}\gamma|x|^{2(\gamma-1)}x_1|\nabla_yu|^2dz
  \\=&\int_{\partial D}F(z,u)\nu_x^1dS-\int_D\frac{\partial F(z,u)}{\partial x_1}dz.
\end{aligned}
\end{equation}
More general, for $1\leq i\leq N$, we have
\begin{equation}\label{x50}
\begin{aligned}
  \frac{1}{p}\int_{\partial D}|\nabla_\gamma u|^p\nu_x^idS
  -&\int_{\partial D}|\nabla_\gamma u|^{p-2}\frac{\partial u}{\partial x_i}\langle\nabla_\gamma u,\nu_\gamma\rangle dS
  -\int_D|\nabla_\gamma u|^{p-2}\gamma|x|^{2(\gamma-1)}x_i|\nabla_yu|^2dz
  \\=&\int_{\partial D}F(z,u)\nu_x^idS-\int_D\frac{\partial F(z,u)}{\partial x_i}dz.
\end{aligned}
\end{equation}

\vskip 0.2cm

\textbf{Similarly, we consider the case of $j=N+1$}.

\vskip 0.2cm

\noindent\textbf{Step 1: Calculate} $\displaystyle\frac{d}{dt}~\mbox{det}~(Jac\Phi^{-1}_t(z))\Big|_{t=0}$.

\vskip 0.2cm

\begin{equation}\label{w1}
\begin{aligned}
  \Phi_t(z)&=z+t\phi(z)e_{N+1}
  \\&=(x_1,\cdots,x_N,y_1+t\phi(z),y_2,\cdots, y_l).
\end{aligned}
\end{equation}
Let
\begin{equation}\label{w2}
  \widetilde{z}=\Phi_t(z),
\end{equation}
i.e.
\begin{equation}\label{w3}
  \widetilde{x_1}=x_1,\;\cdots,\;\widetilde{x_N}=x_N,\;\widetilde{y_1}=y_1+t\phi(z),\;\widetilde{y_2}=y_2\;\cdots,\; \widetilde{y_l}=y_l.
\end{equation}
The Jacobian matrix of $\Phi_t(z)$ is as follows:
\begin{equation}\label{w4}
Jac\Phi_t(z)=
\begin{bmatrix}
1&\cdots&0&0&0&\cdots&0\\
\vdots&\ddots&\vdots&\vdots&\vdots&\ddots&\vdots\\
0&\cdots&1&0&0&\cdots&0\\
t\frac{\partial\phi(z)}{\partial x_1}&\cdots&t\frac{\partial\phi(z)}{\partial x_N}&1+t\frac{\partial\phi(z)}{\partial y_1}&t\frac{\partial\phi(z)}{\partial y_2}&\cdots&t\frac{\partial\phi(z)}{\partial y_l}\\
0&\cdots&0&0&1&\cdots&0\\
\vdots&\ddots&\vdots&\vdots&\vdots&\ddots&\vdots\\
0&\cdots&0&0&0&\cdots&1\\
\end{bmatrix}.
\end{equation}
Then we have
\begin{equation}\label{w5}
  ~\mbox{det}~Jac\Phi_t(z)=1+t\frac{\partial\phi(z)}{\partial y_1}.
\end{equation}
By relationship \eqref{x30} and differentiate at $t=0$ we obtain
\begin{equation}\label{w6}
  \frac{d}{dt}~\mbox{det}~(Jac\Phi^{-1}_t(z))\Big|_{t=0}=-\frac{\partial\phi(z)}{\partial y_1}.
\end{equation}

\vskip 0.2cm

\noindent\textbf{Step 2: Calculate $\displaystyle\frac{1}{p}\frac{d}{dt}\int_\Omega|\nabla_\gamma u(\Phi_t(z))|^pdz\Big|_{t=0}$.}

\vskip 0.2cm

Note that \eqref{x35} holds.
According to \eqref{w1}-\eqref{w3}, by calculation we have
\begin{equation}\label{w7}
  |\nabla_xu(\widetilde{z})|^2=\sum_{i=1}^N\left(\frac{\partial u(\widetilde{z})}{\partial\widetilde{x_i}}
  +\frac{\partial u(\widetilde{z})}{\partial \widetilde{y_1}}t\frac{\partial\phi(z)}{\partial x_i}\right)^2,
\end{equation}

\begin{equation}\label{w8}
  |\nabla_yu(\widetilde{z})|^2=\sum_{j=1}^l\left(\frac{\partial u(\widetilde{z})}{\partial\widetilde{y_j}}
  +\frac{\partial u(\widetilde{z})}{\partial \widetilde{y_1}}t\frac{\partial\phi(z)}{\partial y_j}\right)^2,
\end{equation}
and
\begin{equation}\label{w9}
  |x|^{2\gamma}=|\widetilde{x}|^{2\gamma}.
\end{equation}
Thus it holds
\begin{equation}\label{w10}
  |\nabla_\gamma u(\widetilde{z})|^2=\sum_{i=1}^N\left(\frac{\partial u(\widetilde{z})}{\partial\widetilde{x_i}}
  +\frac{\partial u(\widetilde{z})}{\partial \widetilde{y_1}}t\frac{\partial\phi(z)}{\partial x_i}\right)^2
  +|\widetilde{x}|^{2\gamma}
  \sum_{j=1}^l\left(\frac{\partial u(\widetilde{z})}{\partial\widetilde{y_j}}
  +\frac{\partial u(\widetilde{z})}{\partial \widetilde{y_1}}t\frac{\partial\phi(z)}{\partial y_j}\right)^2,
\end{equation}
and
\begin{equation}\label{w11}
  |\nabla_\gamma u(\widetilde{z})|^p=\left(\sum_{i=1}^N\left(\frac{\partial u(\widetilde{z})}{\partial\widetilde{x_i}}
  +\frac{\partial u(\widetilde{z})}{\partial \widetilde{y_1}}t\frac{\partial\phi(z)}{\partial x_i}\right)^2
  +|\widetilde{x}|^{2\gamma}
  \sum_{j=1}^l\left(\frac{\partial u(\widetilde{z})}{\partial\widetilde{y_j}}
  +\frac{\partial u(\widetilde{z})}{\partial \widetilde{y_1}}t\frac{\partial\phi(z)}{\partial y_j}\right)^2\right)^\frac{p}{2}.
\end{equation}

Accordingly, we have
\begin{equation}\label{w12}
\begin{aligned}
  &\int_\Omega|\nabla_\gamma u(\Phi_t(z))|^pdz
\\=&\int_\Omega\left(\sum_{i=1}^N\left(\frac{\partial u(\widetilde{z})}{\partial\widetilde{x_i}}
  +\frac{\partial u(\widetilde{z})}{\partial \widetilde{y_1}}t\frac{\partial\phi(z)}{\partial x_i}\right)^2
  +|\widetilde{x}|^{2\gamma}
  \sum_{j=1}^l\left(\frac{\partial u(\widetilde{z})}{\partial\widetilde{y_j}}
  +\frac{\partial u(\widetilde{z})}{\partial \widetilde{y_1}}t\frac{\partial\phi(z)}{\partial y_j}\right)^2\right)^\frac{p}{2}dz
\\=&\int_\Omega\left(\sum_{i=1}^N\left(\frac{\partial u(\widetilde{z})}{\partial\widetilde{x_i}}
  +\frac{\partial u(\widetilde{z})}{\partial \widetilde{y_1}}t\frac{\partial\phi(\Phi_t^{-1}(\widetilde{z}))}{\partial (\Phi_t^{-1}(\widetilde{z}))_i}\right)^2
  +|\widetilde{x}|^{2\gamma}
  \sum_{j=1}^l\left(\frac{\partial u(\widetilde{z})}{\partial\widetilde{y_j}}
  +\frac{\partial u(\widetilde{z})}{\partial \widetilde{y_1}}t\frac{\partial\phi(\Phi_t^{-1}(\widetilde{z}))}{\partial (\Phi_t^{-1}(\widetilde{z}))_{N+j}}\right)^2\right)^\frac{p}{2}
  \\&~\mbox{det}~(Jac\Phi^{-1}_t(\widetilde{z}))d\widetilde{z}
\\=&\int_\Omega\left(\sum_{i=1}^N\left(\frac{\partial u(z)}{\partial x_i}
  +\frac{\partial u(z)}{\partial y_1}t\frac{\partial\phi(\Phi_t^{-1}(z))}{\partial (\Phi_t^{-1}(z))_i}\right)^2
  +|x|^{2\gamma}
  \sum_{j=1}^l\left(\frac{\partial u(z)}{\partial y_j}
  +\frac{\partial u(z)}{\partial y_1}t\frac{\partial\phi(\Phi_t^{-1}(z))}{\partial (\Phi_t^{-1}(z))_{N+j}}\right)^2\right)^\frac{p}{2}
  \\&~\mbox{det}~(Jac\Phi^{-1}_t(z))dz.
\end{aligned}
\end{equation}

Then by \eqref{x35}, we have
\begin{equation}\label{w13}
\begin{aligned}
  &\frac{d}{dt}\int_\Omega|\nabla_\gamma u(\Phi_t(z))|^pdz
\\=&\int_\Omega\frac{d}{dt}~\mbox{det}~(Jac\Phi^{-1}_t(z))
\\\cdot&\left(\sum_{i=1}^N\left(\frac{\partial u(z)}{\partial x_i}
  +\frac{\partial u(z)}{\partial y_1}t\frac{\partial\phi(\Phi_t^{-1}(z))}{\partial (\Phi_t^{-1}(z))_i}\right)^2
  +|x|^{2\gamma}
  \sum_{j=1}^l\left(\frac{\partial u(z)}{\partial y_j}
  +\frac{\partial u(z)}{\partial y_1}t\frac{\partial\phi(\Phi_t^{-1}(z))}{\partial (\Phi_t^{-1}(z))_{N+j}}\right)^2\right)^\frac{p}{2}
\\+&~\mbox{det}~(Jac\Phi^{-1}_t(z))
\\\cdot&\frac{p}{2}\left(\sum_{i=1}^N\left(\frac{\partial u(z)}{\partial x_i}
  +\frac{\partial u(z)}{\partial y_1}t\frac{\partial\phi(\Phi_t^{-1}(z))}{\partial (\Phi_t^{-1}(z))_i}\right)^2
  +|x|^{2\gamma}
  \sum_{j=1}^l\left(\frac{\partial u(z)}{\partial y_j}
  +\frac{\partial u(z)}{\partial y_1}t\frac{\partial\phi(\Phi_t^{-1}(z))}{\partial (\Phi_t^{-1}(z))_{N+j}}\right)^2\right)^{\frac{p}{2}-1}
\\\cdot&\left(\sum_{i=1}^N2\left(\frac{\partial u(z)}{\partial x_i}
  +\frac{\partial u(z)}{\partial y_1}t\frac{\partial\phi(\Phi_t^{-1}(z))}{\partial (\Phi_t^{-1}(z))_i}\right)\frac{\partial u(z)}{\partial y_1}\frac{\partial\phi(\Phi_t^{-1}(z))}{\partial (\Phi_t^{-1}(z))_i}\right.
\\+&\left.|x|^{2\gamma}\sum_{j=1}^l2\left(\frac{\partial u(z)}{\partial y_j}
  +\frac{\partial u(z)}{\partial y_1}t\frac{\partial\phi(\Phi_t^{-1}(z))}{\partial(\Phi_t^{-1}(z))_{N+j}}\right)  \frac{\partial u(z)}{\partial y_1}\frac{\partial\phi(\Phi_t^{-1}(z))}{\partial(\Phi_t^{-1}(z))_{N+j}}\right)dz.
\end{aligned}
\end{equation}

Since $\Phi_0(z)=z$, then $\Phi_0^{-1}(z)=z$ and $~\mbox{det}~Jac\Phi_0^{-1}(z)=1$. And by \textbf{Step 1}, we get
\begin{equation}\label{w14}
\begin{aligned}
  &\frac{d}{dt}\int_\Omega|\nabla_\gamma u(\Phi_t(z))|^pdz\Big|_{t=0}
  \\=&\int_\Omega-\frac{\partial\phi(z)}{\partial y_1}\left(\sum_{i=1}^N\left(\frac{\partial u(z)}{\partial x_i}\right)^2+|x|^{2\gamma}\sum_{j=l}^l\left(\frac{\partial u(z)}{\partial y_j}\right)^2\right)^{\frac{p}{2}}
  \\+&\frac{p}{2}\left(\sum_{i=1}^N\left(\frac{\partial u(z)}{\partial x_i}\right)^2+|x|^{2\gamma}\sum_{j=l}^l\left(\frac{\partial u(z)}{\partial y_j}\right)^2\right)^{\frac{p}{2}-1}
\\\cdot&\left(\sum_{i=1}^N2\frac{\partial u(z)}{\partial x_i}
  \frac{\partial u(z)}{\partial y_1}\frac{\partial\phi(z)}{\partial x_i}
+|x|^{2\gamma}\sum_{j=1}^l2\frac{\partial u(z)}{\partial y_j}
 \frac{\partial u(z)}{\partial y_1}\frac{\partial\phi(z)}{\partial y_j}\right)dz
 \\=&\int_\Omega-\frac{\partial\phi(z)}{\partial y_1}|\nabla_\gamma u|^p+p|\nabla_\gamma u|^{p-2}\frac{\partial u}{\partial y_1}
 \langle\nabla_\gamma u,\nabla_\gamma\phi\rangle dz.
\end{aligned}
\end{equation}
Similar to \eqref{x37} and \eqref{x39}, we obtain
\begin{equation}\label{w16}
  \int_\Omega-\frac{\partial\phi(z)}{\partial y_1}|\nabla_\gamma u|^pdz
 \overset{\delta\rightarrow0}\longrightarrow\int_{\partial D}|\nabla_\gamma u|^p\nu_y^1dS,
\end{equation}
and
\begin{equation}\label{w17}
\begin{aligned}
  &\int_\Omega|\nabla_\gamma u|^{p-2}\frac{\partial u}{\partial y_1}\langle\nabla_\gamma u,\nabla_\gamma \phi\rangle dz
\overset{\delta\rightarrow0}\longrightarrow-\int_{\partial D}|\nabla_\gamma u|^{p-2}\frac{\partial u}{\partial y_1}\langle\nabla_\gamma u,\nu_\gamma\rangle dS.
\end{aligned}
\end{equation}
Therefore, by \eqref{w14}-\eqref{w17} we have
\begin{equation}\label{w18}
\begin{aligned}
  &\frac{1}{p}\frac{d}{dt}\int_\Omega|\nabla_\gamma u(\Phi_t(z))|^pdz\Big|_{t=0}
  \\\overset{\delta\rightarrow0}\longrightarrow&\frac{1}{p}\int_{\partial D}|\nabla_\gamma u|^p\nu_y^1dS
  -\int_{\partial D}|\nabla_\gamma u|^{p-2}\frac{\partial u}{\partial y_1}\langle\nabla_\gamma u,\nu_\gamma\rangle dS.
\end{aligned}
\end{equation}

\vskip 0.2cm

\noindent\textbf{Step 3: Calculate $\displaystyle\frac{d}{dt}\int_\Omega f(z,u(\Phi_t(z)))u(\Phi_t(z))dz\Big|_{t=0}$.}

\vskip 0.2cm

Similar to \textbf{Step 3} in case of $j=1$, we have
\begin{equation}\label{w19}
  \frac{d}{dt}\int_\Omega f(z,u(\Phi_t(z)))u(\Phi_t(z))dz\Big|_{t=0}\overset{\delta\rightarrow0}\longrightarrow\int_{\partial D}F(z,u)\nu_y^1dS-\int_D\frac{\partial F(z,u)}{\partial y_1}dz.
\end{equation}

\vskip 0.2cm

\noindent\textbf{Step 4: Complete the proof.}

\vskip 0.2cm

According to \eqref{x17}, \eqref{w18} and \eqref{w19}, we have
\begin{equation}\label{w20}
\begin{aligned}
  \frac{1}{p}\int_{\partial D}|\nabla_\gamma u|^p\nu_y^1dS
  -&\int_{\partial D}|\nabla_\gamma u|^{p-2}\frac{\partial u}{\partial y_1}\langle\nabla_\gamma u,\nu_\gamma\rangle dS
=\int_{\partial D}F(z,u)\nu_y^1dS-\int_D\frac{\partial F(z,u)}{\partial y_1}dz.
\end{aligned}
\end{equation}
More general, for $1\leq j\leq l$, we have
\begin{equation}\label{w21}
\begin{aligned}
  \frac{1}{p}\int_{\partial D}|\nabla_\gamma u|^p\nu_y^jdS
  -&\int_{\partial D}|\nabla_\gamma u|^{p-2}\frac{\partial u}{\partial y_j}\langle\nabla_\gamma u,\nu_\gamma\rangle dS
=\int_{\partial D}F(z,u)\nu_y^jdS-\int_D\frac{\partial F(z,u)}{\partial y_j}dz.
\end{aligned}
\end{equation}

\end{proof}

\begin{Rem}
Note that
\begin{equation}\label{note1}
  \langle\nabla_\gamma u,\nu_\gamma\rangle=\frac{\partial u}{\partial \nu_x}+|x|^{2\gamma}\frac{\partial u}{\partial \nu_y}.
\end{equation}
When $p=2$, our result is the same as that of the proposition from \cite{xiaodong}.
\end{Rem}

\begin{Prop}\label{Prop1}\textbf{(Theorem 1.1, \cite{xiaodong})}\\
If $u\in H^{1,2}_\gamma(\Omega)\cap C^2(\Omega)$ is the solution of
\begin{equation}\label{i1}
  -\Delta_\gamma u=f(z,u),\;\;~\mbox{in}~\Omega\subset\R^{N+l},
\end{equation}
then $u$ satisfies
\begin{equation}\label{.ma2}
\begin{aligned}
\frac{1}{2}\int_{\partial D}|\nabla_\gamma u|^2\nu_x^idS-&\int_{\partial D}\left(\frac{\partial u}{\partial \nu_x}+|x|^{2\gamma}\frac{\partial u}{\partial \nu_y}\right)\frac{\partial u}{\partial x_i}dS-\gamma\int_D|\nabla_y u|^2|x|^{2(\gamma-1)}x_idz
\\=&\int_{\partial D}F(z,u)\nu_x^idS-\int_D\frac{\partial F(z,u)}{\partial x_i}dz,\;\;\;\;~\mbox{i=1,\ldots,N},
\end{aligned}
\end{equation}
and
\begin{equation}\label{.ma3}
\begin{aligned}
\\&\frac{1}{2}\int_{\partial D}|\nabla_\gamma u|^2\nu_y^jdS-\int_{\partial D}\left(\frac{\partial u}{\partial \nu_x}+|x|^{2\gamma}\frac{\partial u}{\partial \nu_y}\right)\frac{\partial u}{\partial y_j}dS
\\=&\int_{\partial D}F(z,u)\nu_y^jdS-\int_D\frac{\partial F(z,u)}{\partial y_j}dz,\;\;\;\;~\mbox{j=1,\ldots,l},
\end{aligned}
\end{equation}
where $D\subset\Omega$ is any domain of $\Omega$, $\nu=(\nu_x,\nu_y)$ is the unit outward normal of the point of $\partial D$, $F(z,u)=\displaystyle\int_0^uf(z,s)ds$, and $F_{x_i}(z,u)=\displaystyle\frac{\partial F(z,u)}{\partial x_i}+f(z,u)\frac{\partial u}{\partial x_i}$.

\end{Prop}

\section{Local Pohozaev identity of scaling type}

In this section, we prove Pohozaev identity generated from scaling of \eqref{x6} in Theorem \ref{th2}.

\begin{proof}[\textbf{Proof of Theorem \ref{th2}:}]

Suppose $D\subset\subset\Omega$, and for any $\delta>0$ small enough define $D_\delta$ as \eqref{x14}.
Now, we take
\begin{equation}\label{d1}
\begin{aligned}
  \Phi_t(z)&=(1+t\phi(z))z
 \\&=((1+t\phi(z))x_1,\cdots,(1+t\phi(z))x_N,(1+t\phi(z))y_1,\cdots,(1+t\phi(z))y_l),
\end{aligned}
\end{equation}
where $t\in\R$, $z=(x,y)=(x_1,\cdots,x_N,y_1,\cdots,y_l)\in\Omega$, $\phi\in W^{1,p}_{\lambda,0}(\Omega)$ satisfies \eqref{x13}, \eqref{x15}, \eqref{x38}. Here, $\nu=(\nu_x,\nu_y)$ is the outward unit normal of $\partial D$ at $\overline{z}=(\overline{x},\overline{y})\in\partial D$ in the sense of \eqref{x16}. Clearly, $\Phi_t(z):\Omega\rightarrow\Omega$ satisfies $\Phi_0(z)=z$ and $\Phi_t(z):\partial\Omega\rightarrow\partial\Omega$.

By the definition \eqref{x10} of functional $I$, we have \eqref{x180}. And it follows from \eqref{x11} that we have \eqref{x17}.

\vskip 0.2cm

\noindent\textbf{Step 1: Calculate} $\displaystyle\frac{d}{dt}~\mbox{det}~(Jac\Phi^{-1}_t(z))\Big|_{t=0}$.

\vskip 0.2cm

Let
\begin{equation}\label{d2}
  \widetilde{z}=\Phi_t(z),
\end{equation}
i.e.
\begin{equation}\label{d3}
  \widetilde{x_1}=(1+t\phi(z))x_1,\;\cdots,\;\widetilde{x_N}=(1+t\phi(z))x_N,
  \;\widetilde{y_1}=(1+t\phi(z))y_1,\;\cdots,\; \widetilde{y_l}=(1+t\phi(z))y_l.
\end{equation}
The Jacobian matrix of $\Phi_t(z)$ is as follows:
\begin{equation}\label{d4}
Jac\Phi_t(z)=
\begin{bmatrix}
1+t\phi(z)+x_1t\frac{\partial\phi(z)}{\partial x_1} & x_1t\frac{\partial\phi(z)}{\partial x_2} & \cdots & x_1t\frac{\partial\phi(z)}{\partial y_l}\\
x_2t\frac{\partial\phi(z)}{\partial x_1} & 1+t\phi(z)+x_2t\frac{\partial\phi(z)}{\partial x_2} & \cdots & x_2t\frac{\partial\phi(z)}{\partial y_l}\\
\vdots & \vdots & \ddots & \vdots\\
y_lt\frac{\partial\phi}{\partial x_1} & y_lt\frac{\partial\phi}{\partial x_2} & \cdots & 1+t\phi(z)+y_lt\frac{\partial\phi(z)}{\partial y_l}\\
\end{bmatrix}.
\end{equation}
From \eqref{d4}, we derive that
\begin{equation}\label{d5}
  ~\mbox{det}~Jac\Phi_t(z)=(1+t\phi(z))^{N+l}+(1+t\phi(z))^{N+l-1}t\langle z,\nabla\phi\rangle.
\end{equation}
According to relation \eqref{x30}, then we have
\begin{equation}\label{d6}
  ~\mbox{det}~Jac\Phi_t^{-1}(z)=\frac{1}{(1+t\phi(z))^{N+l}+(1+t\phi(z))^{N+l-1}t\langle z,\nabla\phi\rangle}.
\end{equation}
Taking the derivative of equation \eqref{d6} with respect to $t$, we get
\begin{equation}\label{d7}
\begin{aligned}
  &\frac{d}{dt}~\mbox{det}~Jac\Phi_t^{-1}(z)
  \\=&-\frac
  {(N+l)(1+t\phi(z))^{N+l-1}\phi(z)}
  {((1+t\phi(z))^{N+l}+(1+t\phi(z))^{N+l-1}t\langle z,\nabla\phi\rangle)^2}
    \\&-\frac
  {(N+l-1)(1+t\phi(z))^{N+l-2}\phi(z)t\langle z,\nabla\phi\rangle
  +(1+t\phi(z))^{N+l-1}\langle z,\nabla\phi\rangle}
  {((1+t\phi(z))^{N+l}+(1+t\phi(z))^{N+l-1}t\langle z,\nabla\phi\rangle)^2}.
\end{aligned}
\end{equation}
Therefore, we obtain
\begin{equation}\label{d8}
  \frac{d}{dt}~\mbox{det}~(Jac\Phi^{-1}_t(z))\Big|_{t=0}=-\left((N+l)\phi(z)+\langle z,\nabla\phi(z)\rangle\right).
\end{equation}

\vskip 0.2cm

\noindent\textbf{Step 2: Calculate $\displaystyle\frac{1}{p}\frac{d}{dt}\int_\Omega|\nabla_\gamma u(\Phi_t(z))|^pdz\Big|_{t=0}$.}

\vskip 0.2cm

Note that it holds \eqref{x35}.
According to \eqref{d1}-\eqref{d3}, by calculation we have
\begin{equation}\label{d9}
  |\nabla_xu(\widetilde{z})|^2=\sum_{i=1}^N\left(\frac{\partial u(\widetilde{z})}{\partial\widetilde{x_i}}
  (1+t\phi(z))+t\frac{\partial\phi(z)}{\partial x_i}\langle\nabla_{\widetilde{z}} u(\widetilde{z}),z\rangle\right)^2,
\end{equation}

\begin{equation}\label{d10}
  |\nabla_yu(\widetilde{z})|^2=\sum_{j=1}^l\left(\frac{\partial u(\widetilde{z})}{\partial\widetilde{y_j}}
  (1+t\phi(z))+t\frac{\partial\phi(z)}{\partial y_j}\langle\nabla_{\widetilde{z}}u(\widetilde{z}),z\rangle\right)^2,
\end{equation}
and
\begin{equation}\label{d11}
  |x|^{2\gamma}=\frac{|\widetilde{x}|^{2\gamma}}{(1+t\phi(z))^{2\gamma}}.
\end{equation}
Thus, we derive that
\begin{equation}\label{d11}
\begin{aligned}
  |\nabla_\gamma u(\widetilde{z})|^2=&\sum_{i=1}^N\left(\frac{\partial u(\widetilde{z})}{\partial\widetilde{x_i}}
  (1+t\phi(z))+t\frac{\partial\phi(z)}{\partial x_i}\langle\nabla_{\widetilde{z}} u(\widetilde{z}),z\rangle\right)^2
  \\+&\frac{|\widetilde{x}|^{2\gamma}}{(1+t\phi(z))^{2\gamma}}
  \sum_{j=1}^l\left(\frac{\partial u(\widetilde{z})}{\partial\widetilde{y_j}}
  (1+t\phi(z))+t\frac{\partial\phi(z)}{\partial y_j}\langle\nabla_{\widetilde{z}}u(\widetilde{z}),z\rangle\right)^2,
\end{aligned}
\end{equation}
and
\begin{equation}\label{d12}
\begin{aligned}
  |\nabla_\gamma u(\widetilde{z})|^p&=\left(\sum_{i=1}^N\left(\frac{\partial u(\widetilde{z})}{\partial\widetilde{x_i}}
  (1+t\phi(z))+t\frac{\partial\phi(z)}{\partial x_i}\langle\nabla_{\widetilde{z}} u(\widetilde{z}),z\rangle\right)^2\right.
  \\&\left.+\frac{|\widetilde{x}|^{2\gamma}}{(1+t\phi(z))^{2\gamma}}
  \sum_{j=1}^l\left(\frac{\partial u(\widetilde{z})}{\partial\widetilde{y_j}}
  (1+t\phi(z))+t\frac{\partial\phi(z)}{\partial y_j}\langle\nabla_{\widetilde{z}}u(\widetilde{z}),z\rangle\right)^2\right)^{\frac{p}{2}}.
\end{aligned}
\end{equation}
Accordingly, we have
\begin{equation}\label{d13}
\begin{aligned}
  &\int_\Omega|\nabla_\gamma u(\Phi_t(z))|^pdz
  \\=&\int_\Omega\left(\sum_{i=1}^N\left(\frac{\partial u(\widetilde{z})}{\partial\widetilde{x_i}}
  (1+t\phi(z))+t\frac{\partial\phi(z)}{\partial x_i}\langle\nabla_{\widetilde{z}} u(\widetilde{z}),z\rangle\right)^2\right.
  \\&\left.+\frac{|\widetilde{x}|^{2\gamma}}{(1+t\phi(z))^{2\gamma}}
  \sum_{j=1}^l\left(\frac{\partial u(\widetilde{z})}{\partial\widetilde{y_j}}
  (1+t\phi(z))+t\frac{\partial\phi(z)}{\partial y_j}\langle\nabla_{\widetilde{z}}u(\widetilde{z}),z\rangle\right)^2\right)^{\frac{p}{2}}dz
  \\=&\int_\Omega\left(\sum_{i=1}^N\left(\frac{\partial u(\widetilde{z})}{\partial\widetilde{x_i}}
  (1+t\phi(\Phi_t^{-1}(\widetilde{z})))+t\frac{\partial\phi(\Phi_t^{-1}(\widetilde{z}))}{\partial (\Phi_t^{-1}(\widetilde{z}))_i}\langle\nabla_{\widetilde{z}} u(\widetilde{z}),\Phi_t^{-1}(\widetilde{z})\rangle\right)^2\right.
  \\&\left.+\frac{|\widetilde{x}|^{2\gamma}}{(1+t\phi(\Phi_t^{-1}(\widetilde{z})))^{2\gamma}}
  \sum_{j=1}^l\left(\frac{\partial u(\widetilde{z})}{\partial\widetilde{y_j}}
  (1+t\phi(\Phi_t^{-1}(\widetilde{z})))+t\frac{\partial\phi(\Phi_t^{-1}(\widetilde{z}))}{\partial (\Phi_t^{-1}(\widetilde{z}))_{N+j}}\langle\nabla_{\widetilde{z}}u(\widetilde{z}),\Phi_t^{-1}(\widetilde{z})
  \rangle\right)^2\right)^{\frac{p}{2}}
  \\&~\mbox{det}~(Jac\Phi^{-1}_t(\widetilde{z}))d\widetilde{z}
  \\=&\int_\Omega\left(\sum_{i=1}^N\left(\frac{\partial u(z)}{\partial x_i}
  (1+t\phi(\Phi_t^{-1}(z)))+t\frac{\partial\phi(\Phi_t^{-1}(z))}{\partial (\Phi_t^{-1}(z))_i}\langle\nabla u(z),\Phi_t^{-1}(z)\rangle\right)^2\right.
  \\&\left.+\frac{|x|^{2\gamma}}{(1+t\phi(\Phi_t^{-1}(z)))^{2\gamma}}
  \sum_{j=1}^l\left(\frac{\partial u(z)}{\partial y_j}
  (1+t\phi(\Phi_t^{-1}(z)))+t\frac{\partial\phi(\Phi_t^{-1}(z))}{\partial (\Phi_t^{-1}(z))_{N+j}}\langle\nabla u(z),\Phi_t^{-1}(z)
  \rangle\right)^2\right)^{\frac{p}{2}}
  \\&~\mbox{det}~(Jac\Phi^{-1}_t(z))dz.
\end{aligned}
\end{equation}

Then by \eqref{x35}, we have
\begin{equation*}
\begin{aligned}
 &\frac{d}{dt}\int_\Omega|\nabla_\gamma u(\Phi_t(z))|^pdz
\\=&\int_\Omega\frac{d}{dt}~\mbox{det}~(Jac\Phi^{-1}_t(z))
\cdot\left(\sum_{i=1}^N\left(\frac{\partial u(z)}{\partial x_i}
  (1+t\phi(\Phi_t^{-1}(z)))+t\frac{\partial\phi(\Phi_t^{-1}(z))}{\partial (\Phi_t^{-1}(z))_i}\langle\nabla u(z),\Phi_t^{-1}(z)\rangle\right)^2\right.
  \\&\left.+\frac{|x|^{2\gamma}}{(1+t\phi(\Phi_t^{-1}(z)))^{2\gamma}}
  \sum_{j=1}^l\left(\frac{\partial u(z)}{\partial y_j}
  (1+t\phi(\Phi_t^{-1}(z)))+t\frac{\partial\phi(\Phi_t^{-1}(z))}{\partial (\Phi_t^{-1}(z))_{N+j}}\langle\nabla u(z),\Phi_t^{-1}(z)
  \rangle\right)^2\right)^{\frac{p}{2}}dz
\end{aligned}
\end{equation*}

\begin{equation}\label{d14}
\begin{aligned}
+&\int_\Omega~\mbox{det}~(Jac\Phi^{-1}_t(z))
\cdot\frac{p}{2}\left(\sum_{i=1}^N\left(\frac{\partial u(z)}{\partial x_i}
  (1+t\phi(\Phi_t^{-1}(z)))+t\frac{\partial\phi(\Phi_t^{-1}(z))}{\partial (\Phi_t^{-1}(z))_i}\langle\nabla u(z),\Phi_t^{-1}(z)\rangle\right)^2\right.
  \\&\left.+\frac{|x|^{2\gamma}}{(1+t\phi(\Phi_t^{-1}(z)))^{2\gamma}}
  \sum_{j=1}^l\left(\frac{\partial u(z)}{\partial y_j}
  (1+t\phi(\Phi_t^{-1}(z)))+t\frac{\partial\phi(\Phi_t^{-1}(z))}{\partial (\Phi_t^{-1}(z))_{N+j}}\langle\nabla u(z),\Phi_t^{-1}(z)
  \rangle\right)^2\right)^{\frac{p}{2}-1}
\\\cdot&\left(\sum_{i=1}^N2\left(\frac{\partial u(z)}{\partial x_i}
  (1+t\phi(\Phi_t^{-1}(z)))+t\frac{\partial\phi(\Phi_t^{-1}(z))}{\partial (\Phi_t^{-1}(z))_i}\langle\nabla u(z),\Phi_t^{-1}(z)\rangle\right)\right.
 \\&\left.\left(\frac{\partial u(z)}{\partial x_i}\phi(\Phi_t^{-1}(z))+\frac{\partial\phi(\Phi_t^{-1}(z))}{\partial (\Phi_t^{-1}(z))_i}\langle\nabla u(z),\Phi_t^{-1}(z)\rangle\right)\right.
    +\left.\frac{-2\gamma|x|^{2\gamma}\phi(\Phi_t^{-1}(z))}{(1+t\phi(\Phi_t^{-1}(z)))^{2\gamma+1}}\right.
    \\&\cdot\left.\sum_{j=1}^l\left(\frac{\partial u(z)}{\partial y_j}
  (1+t\phi(\Phi_t^{-1}(z)))+t\frac{\partial\phi(\Phi_t^{-1}(z))}{\partial (\Phi_t^{-1}(z))_{N+j}}\langle\nabla u(z),\Phi_t^{-1}(z)
  \rangle\right)^2\right.
    \\&+\left.\frac{|x|^{2\gamma}}{(1+t\phi(\Phi_t^{-1}(z)))^{2\gamma}}\right.
    \cdot\left.\sum_{j=1}^l2\left(\frac{\partial u(z)}{\partial y_j}
  (1+t\phi(\Phi_t^{-1}(z)))+t\frac{\partial\phi(\Phi_t^{-1}(z))}{\partial (\Phi_t^{-1}(z))_{N+j}}\langle\nabla u(z),\Phi_t^{-1}(z)
  \rangle\right)\right.
  \\&\left.\left(\frac{\partial u(z)}{\partial y_j}
  \phi(\Phi_t^{-1}(z))+\frac{\partial\phi(\Phi_t^{-1}(z))}{\partial (\Phi_t^{-1}(z))_{N+j}}\langle\nabla u(z),\Phi_t^{-1}(z)
  \rangle\right)\right)dz.
\end{aligned}
\end{equation}

Since $\Phi_0(z)=z$, then $\Phi_0^{-1}(z)=z$ and $~\mbox{det}~Jac\Phi_0^{-1}(z)=1$. And by \textbf{Step 1}, we get

\begin{equation}\label{d15}
\begin{aligned}
 &\frac{d}{dt}\int_\Omega|\nabla_\gamma u(\Phi_t(z))|^pdz\Big|_{t=0}
\\=&-\left((N+l)\phi(z)+\langle z,\nabla\phi(z)\rangle\right)
\left(\sum_{i=1}^N\left(\frac{\partial u(z)}{\partial x_i}\right)^2+|x|^{2\gamma}
  \sum_{j=1}^l\left(\frac{\partial u(z)}{\partial y_j}\right)^2\right)^{\frac{p}{2}}
            \\&+\frac{p}{2}\left(\sum_{i=1}^N\left(\frac{\partial u(z)}{\partial x_i}\right)^2+|x|^{2\gamma}
            \sum_{j=1}^l\left(\frac{\partial u(z)}{\partial y_j}\right)^2\right)^{\frac{p}{2}-1}
\cdot\left(\sum_{i=1}^N2\frac{\partial u(z)}{\partial x_i}
\left(\frac{\partial u(z)}{\partial x_i}\phi(z)+\frac{\partial\phi(z)}{\partial x_i}\langle\nabla u,z\rangle\right)\right.
             \\&\left.+(-2\gamma)|x|^{2\gamma}\phi(z)\sum_{j=1}^l\left(\frac{\partial u(z)}{\partial y_j}\right)^2\right.
  \left.+|x|^{2\gamma}
  \sum_{j=1}^l2\frac{\partial u(z)}{\partial y_j}
  \left(\frac{\partial u(z)}{\partial y_j}\phi(z)+\frac{\partial\phi(z)}{\partial y_j}\langle\nabla u,z\rangle\right)\right)
\\=&-\left((N+l)\phi(z)+\langle z,\nabla\phi(z)\rangle\right)|\nabla_\gamma u|^p
\\&+p|\nabla_\gamma u|^{p-2}
\left(\phi(z)|\nabla_\gamma u|^2+\langle\nabla u,z\rangle\langle\nabla_\gamma u,\nabla_\gamma \phi\rangle
-\gamma|x|^{2\gamma}\phi(z)|\nabla_yu|^2\right)
\\=&(-(N+l)+p)|\nabla_\gamma u|^p\phi(z)-|\nabla_\gamma u|^p\langle z,\nabla\phi\rangle
\\&+p|\nabla_\gamma u|^{p-2}\langle\nabla u,z\rangle\langle\nabla_\gamma u,\nabla_\gamma\phi\rangle
-p|\nabla_\gamma u|^{p-2}\gamma|x|^{2\gamma}\phi|\nabla_yu|^2.
\end{aligned}
\end{equation}
In which, by \eqref{x13} and \eqref{x38} we obtain
\begin{equation}\label{d17}
  \int_\Omega|\nabla_\gamma u|^p\phi(z)dz
 =\int_{D_\delta}|\nabla_\gamma u|^p\phi(z)dz
 \overset{\delta\rightarrow0}\longrightarrow\int_D|\nabla_\gamma u|^pdz,
\end{equation}
and
\begin{equation}\label{d18}
\begin{aligned}
  &\int_\Omega|\nabla_\gamma u|^p\langle z,\nabla\phi\rangle dz
\\=&\int_{D_\delta\setminus D}|\nabla_\gamma u|^p\langle z,\nabla\phi\rangle dz
\\\overset{\delta\rightarrow0}\longrightarrow&-\int_{\partial D}|\nabla_\gamma u|^p\langle z,\nu\rangle dS,
\end{aligned}
\end{equation}
and
\begin{equation}\label{d19}
\begin{aligned}
  &\int_\Omega|\nabla_\gamma u|^{p-2}\langle\nabla u,z\rangle\langle\nabla_\gamma u,\nabla_\gamma\phi\rangle dz
  \\=&\int_{D_\delta\setminus D}|\nabla_\gamma u|^{p-2}\langle\nabla u,z\rangle\langle\nabla_\gamma u,\nabla_\gamma\phi\rangle dz
  \\\overset{\delta\rightarrow0}\longrightarrow&-\int_{\partial D}|\nabla_\gamma u|^{p-2}\langle\nabla u,z\rangle\langle\nabla_\gamma u,\nu_\gamma\rangle dS,
\end{aligned}
\end{equation}
and
\begin{equation}\label{d20}
\begin{aligned}
  &\int_\Omega|\nabla_\gamma u|^{p-2}\gamma|x|^{2\gamma}\phi|\nabla_yu|^2dz
  \\=&\int_{D_\delta}|\nabla_\gamma u|^{p-2}\gamma|x|^{2\gamma}\phi|\nabla_yu|^2dz
  \\\overset{\delta\rightarrow0}\longrightarrow&\int_D|\nabla_\gamma u|^{p-2}\gamma|x|^{2\gamma}|\nabla_yu|^2dz.
\end{aligned}
\end{equation}

Therefore, by \eqref{d15}-\eqref{d20} we have
\begin{equation}\label{d21}
\begin{aligned}
  &\frac{1}{p}\frac{d}{dt}\int_\Omega|\nabla_\gamma u(\Phi_t(z))|^pdz\Big|_{t=0}
  \\\overset{\delta\rightarrow0}\longrightarrow&\left(-\frac{N+l}{p}+1\right)\int_D|\nabla_\gamma u|^pdz
  +\frac{1}{p}\int_{\partial D}|\nabla_\gamma u|^p\langle z,\nu\rangle dS
 \\-&\int_{\partial D}|\nabla_\gamma u|^{p-2}\langle\nabla u,z\rangle\langle\nabla_\gamma u,\nu_\gamma\rangle dS
-\int_D|\nabla_\gamma u|^{p-2}\gamma|x|^{2\gamma}|\nabla_yu|^2dz.
\end{aligned}
\end{equation}

\vskip 0.2cm

\noindent\textbf{Step 3: Calculate $\displaystyle\frac{d}{dt}\int_\Omega g(z)u(\Phi_t(z))dz\Big|_{t=0}$.}

\vskip 0.2cm

Note that it holds \eqref{x41}. Since
\begin{equation}\label{d22}
  u(\Phi_t(z))=u((1+t\phi(z))z),
\end{equation}
and by \eqref{d3}, it is easy to find that
\begin{equation}\label{d23}
  \frac{d}{dt}u(\Phi_t(z))=\phi(z)\langle\nabla u(\Phi_t(z)),z\rangle,
\end{equation}
and thus
\begin{equation}\label{d24}
  \frac{d}{dt}u(\Phi(z))\Big|_{t=0}=\phi(z)\langle\nabla u(z),z\rangle.
\end{equation}
Then
\begin{equation}\label{d25}
\begin{aligned}
  &\int_\Omega g(z)\frac{d}{dt}u(\Phi_t(z))\Big|_{t=0}dz
\\=&\int_\Omega g(z)\phi(z)\langle\nabla u,z\rangle dz
\\=&\int_{D_\delta} g(z)\phi(z)\langle\nabla u,z\rangle dz
\\\overset{\delta\rightarrow0}\longrightarrow&\int_D g(z)\langle\nabla u,z\rangle dz.
\end{aligned}
\end{equation}
On the other hand, by the definition of \eqref{x420}, we have
\begin{equation}\label{d26}
\begin{aligned}
  &\int_D g(z)\langle\nabla u,z\rangle dz
  \\=&\int_D f(z,u)\langle\nabla u,z\rangle dz
  \\=&\int_D\langle f(z,u)\nabla u,z\rangle dz
  \\=&\int_D\langle \nabla F(z,u),z\rangle-\langle \nabla_zF(z,u),z\rangle dz
  \\=&\int_{\partial D}F(z,u)\langle z,\nu\rangle dS
  -(N+l)\int_DF(z,u)dz
  -\int_D\langle \nabla_zF(z,u),z\rangle dz,
\end{aligned}
\end{equation}
where $\nabla F(z,u)$ and $\nabla_zF(z,u)$ satisfies \eqref{sjm2}.

Thus from \eqref{d25} and \eqref{d26}, we have
\begin{equation}\label{d27}
\begin{aligned}
  &\int_\Omega g(z)\frac{d}{dt}u(\Phi_t(z))\Big|_{t=0}dz
  \\\overset{\delta\rightarrow0}\longrightarrow&\int_{\partial D}F(z,u)\langle z,\nu\rangle dS
  -(N+l)\int_DF(z,u)dz
  -\int_D\langle \nabla_zF(z,u),z\rangle dz.
\end{aligned}
\end{equation}
Therefore by \eqref{x41} and \eqref{d27}, we have
\begin{equation}\label{d28}
\begin{aligned}
  &\frac{d}{dt}\int_\Omega g(z)u(\Phi_t(z))dz\Big|_{t=0}
   \\\overset{\delta\rightarrow0}\longrightarrow&\int_{\partial D}F(z,u)\langle z,\nu\rangle dS
  -(N+l)\int_DF(z,u)dz
  -\int_D\langle \nabla_zF(z,u),z\rangle dz.
\end{aligned}
\end{equation}

\vskip 0.2cm

\noindent\textbf{Step 4: Complete the proof.}

\vskip 0.2cm

According to \eqref{x17}, \eqref{d21} and \eqref{d28}, we have
\begin{equation}\label{d29}
\begin{aligned}
&\left(-\frac{N+l}{p}+1\right)\int_D|\nabla_\gamma u|^pdz
  +\frac{1}{p}\int_{\partial D}|\nabla_\gamma u|^p\langle z,\nu\rangle dS
 \\&-\int_{\partial D}|\nabla_\gamma u|^{p-2}\langle\nabla u,z\rangle\langle\nabla_\gamma u,\nu_\gamma\rangle dS
-\int_D|\nabla_\gamma u|^{p-2}\gamma|x|^{2\gamma}|\nabla_yu|^2dz
\\=&\int_{\partial D}F(z,u)\langle z,\nu\rangle dS
  -(N+l)\int_DF(z,u)dz
  -\int_D\langle \nabla_zF(z,u),z\rangle dz.
\end{aligned}
\end{equation}

At last, we need to deal with the term $\displaystyle\int_D|\nabla_\gamma u|^pdz$. Since $u\in W_{\lambda,0}^{1,p}(\Omega)$ is a weak solution of \eqref{x6}, that is
\begin{equation}\label{d30}
  \int_\Omega|\nabla_\gamma u|^{p-2}\langle\nabla_\gamma u,\nabla_\gamma\varphi\rangle dz=\int_\Omega f(z,u)\varphi dz,\;\;\forall\varphi\in W_{\lambda,0}^{1,p}(\Omega).
\end{equation}
Take $\varphi=\phi u$ and insert in \eqref{d30}, we obtain
\begin{equation}\label{d31}
   \int_\Omega|\nabla_\gamma u|^{p-2}u\langle\nabla_\gamma u,\nabla_\gamma\phi\rangle dz
   +\int_\Omega|\nabla_\gamma u|^p\phi dz
   =\int_\Omega f(z,u)\phi udz.
\end{equation}
In \eqref{d31}, by \eqref{x13} and \eqref{x38} we have
\begin{equation}\label{d32}
\begin{aligned}
  &\int_\Omega|\nabla_\gamma u|^{p-2}u\langle\nabla_\gamma u,\nabla_\gamma\phi\rangle dz
\\=&\int_{D_\delta\setminus D}|\nabla_\gamma u|^{p-2}u\langle\nabla_\gamma u,\nabla_\gamma\phi\rangle dz
\\\overset{\delta\rightarrow0}\longrightarrow&-\int_{\partial D}|\nabla_\gamma u|^{p-2}u\langle\nabla_\gamma u,\nu_\gamma\rangle dS,
\end{aligned}
\end{equation}
and
\begin{equation}\label{d33}
  \int_\Omega|\nabla_\gamma u|^p\phi dz=\int_{D_\delta}|\nabla_\gamma u|^p\phi dz
  \overset{\delta\rightarrow0}\longrightarrow\int_D|\nabla_\gamma u|^pdz,
\end{equation}
and
\begin{equation}\label{d34}
  \int_\Omega f(z,u)\phi udz=\int_{D_\delta}f(z,u)\phi udz
  \overset{\delta\rightarrow0}\longrightarrow\int_Df(z,u)udz.
\end{equation}
Thus from \eqref{d31}-\eqref{d34}, we know that
\begin{equation}\label{d35}
  \int_D|\nabla_\gamma u|^pdz=\int_Df(z,u)udz+\int_{\partial D}|\nabla_\gamma u|^{p-2}u\langle\nabla_\gamma u,\nu_\gamma\rangle dS.
\end{equation}
Combining \eqref{d29} and \eqref{d35}, finally we have
\begin{equation}\label{d36}
\begin{aligned}
&\left(-\frac{N+l}{p}+1\right)\int_Df(z,u)udz+\left(-\frac{N+l}{p}+1\right)\int_{\partial D}|\nabla_\gamma u|^{p-2}u\langle\nabla_\gamma u,\nu_\gamma\rangle dS
 \\& +\frac{1}{p}\int_{\partial D}|\nabla_\gamma u|^p\langle z,\nu\rangle dS
-\int_{\partial D}|\nabla_\gamma u|^{p-2}\langle\nabla u,z\rangle\langle\nabla_\gamma u,\nu_\gamma\rangle dS
-\int_D|\nabla_\gamma u|^{p-2}\gamma|x|^{2\gamma}|\nabla_yu|^2dz
\\=&\int_{\partial D}F(z,u)\langle z,\nu\rangle dS
  -(N+l)\int_DF(z,u)dz
  -\int_D\langle \nabla_zF(z,u),z\rangle dz.
\end{aligned}
\end{equation}

\end{proof}

\begin{Rem}\label{Rem1}
Note that $\nabla_xu=\frac{\partial u}{\partial\nu}\nu_x$, $\nabla_yu=\frac{\partial u}{\partial\nu}\nu_y$, since $u=0$ on $\partial\Omega$. Thus we have
\begin{equation}\label{yw260}
  \int_{\partial\Omega}\left(\frac{\partial u}{\partial\nu}\right)^2(\nu_x^2+|x|^{2\gamma}\nu_y^2)(x\cdot\nu_x+y\cdot\nu_y)dS=\int_{\partial\Omega}|\nabla_\gamma u|^2(z\cdot\nu)dS.
\end{equation}
Accordingly, we have
\begin{equation}\label{d39}
\begin{aligned}
  &\int_{\partial\Omega}|\nabla_\gamma u|^{p-2}\langle\nabla u,z\rangle\langle\nabla_\gamma u,\nu_\gamma\rangle dS
\\=&\int_{\partial\Omega}|\nabla_\gamma u|^{p-2}
\left(\frac{\partial u}{\partial\nu}\right)^2
(\nu_x+|x|^{2\gamma}\nu_y)(x\cdot\nu_x+y\cdot\nu_y)dS
\\=&\int_{\partial\Omega}|\nabla_\gamma u|^p
\langle z,\nu\rangle dS.
\end{aligned}
\end{equation}
When $p=2$ and take $D=\Omega$, the result we derive is the same as the following proposition from \cite{xiaodong}.
\end{Rem}

\begin{Prop}\label{Prop2}\textbf{(Theorem 1.2, \cite{xiaodong})}\\
If $u\in H^{1,2}_\gamma(\Omega)\cap C^2(\overline{\Omega})$ is the solution of
\begin{equation}\label{ijk4}
\begin{cases}
-\Delta_\gamma u=f(z,u),\;\;~&\mbox{in}~\Omega,\\[2mm]
u=0,\;\;~&\mbox{on}~\partial\Omega,
\end{cases}
\end{equation}
then $u$ satisfies
\begin{equation}\label{.ma7}
\begin{aligned}
  (N+l)\int_\Omega F(z,u)dz&=-\int_\Omega z\cdot\nabla_z F(z,u)dz+\left(\frac{N+l}{2}-1\right)\int_\Omega f(z,u)udz
+\gamma\int_\Omega|x|^{2\gamma}|\nabla_yu|^2dz\\&+\frac{1}{2}\int_{\partial\Omega}\left(\frac{\partial u}{\partial\nu}\right)^2\left(\nu_x^2+|x|^{2\gamma}\nu_y^2\right)\left(x\cdot\nu_x+y\cdot\nu_y\right)dS,
\end{aligned}
\end{equation}
where $\Omega\subset\R^{N+l}$ is any domain of $\R^{N+l}$ and $\nu=(\nu_x,\nu_y)$ is the unit outward normal of the point of $\partial\Omega$, $F(z,u)=\displaystyle\int_0^uf(z,s)ds$ and $\nabla F(z,u)=\nabla_zF(z,u)+f(z,u)\nabla u$.

\end{Prop}

\section{Global Pohozaev identity of scaling type}

In this section, we prove the global Pohozaev identity generated from scaling of \eqref{x6} in Theorem \ref{th3}.

\begin{proof}[\textbf{Proof of Theorem \ref{th3}:}]
For the proof of the first part, one can refer to Remark \ref{Rem1}. By combining \eqref{d36} and \eqref{d39}, we derive \eqref{d37}.

Now we consider $\Omega=\R^{N+l}$. Since $|\nabla_\gamma u|\in L^p(\R^{N+l})$, we have
\begin{equation}\label{zxd1}
\begin{aligned}
  \int_{\R^{N+l}}|\nabla_\gamma u|^pdz
  &=\int_0^{+\infty}\int_{\partial B_r(0)}|\nabla_\gamma u(r,\theta)|^pdSr^{N+l-1}dr
\\&=\int_0^{+\infty}r^{N+l-2}\int_{\partial B_r(0)}|\nabla_\gamma u(r,\theta)|^prdSdr
\\&=\int_0^{+\infty}r^{N+l-2}\int_{\partial B_r(0)}|\nabla_\gamma u(r,\theta)|^p(z\cdot\nu)dSdr<+\infty.
\end{aligned}
\end{equation}
Then there exists a sequence $r_n\rightarrow+\infty$ such that as $n\rightarrow+\infty$,
\begin{equation}\label{zxd2}
  r_n\int_{\partial B_{r_n}(0)}|\nabla_\gamma u(r_n,\theta)|^pdS\rightarrow0.
\end{equation}
Suppose this is false, then we have
\begin{equation}\label{zxd3}
\mathop {\lim \inf }\limits_{r \to  + \infty }r\int_{\partial B_r(0)}|\nabla_\gamma u(r,\theta)|^pdS>\alpha>0.
\end{equation}
According to sign-preserving theorem of limit, $\xi(r)=r\displaystyle\int_{\partial B_r(0)}|\nabla_\gamma u(r,\theta)|^pdS>\alpha>0$, hence
\begin{equation}\label{zxd4}
  \int_0^{+\infty}r^{N+l-2}\xi(r)dr>\alpha\int_0^{+\infty}r^{N+l-2}dr=+\infty,
\end{equation}
which contradicts the fact that $|\nabla_\gamma u|\in L^p(\R^{N+l})$. Thus, combining \eqref{d37} and \eqref{zxd2}, we get \eqref{d38}.

\end{proof}

\textbf{Data Availability Statement}

No data, models, or code were generated or used during this study.
\\*

\textbf{Conflict of interest statement}

The authors declare that they have no conflicts of interest in the research presented in this manuscript.
\\*

\end{document}